\documentclass[12pt,fceqn]{article}
\usepackage{amsmath}
\usepackage{amsfonts}
\usepackage{cite}
\usepackage{amssymb,graphics,graphicx,subfigure,caption2,rotating}
\usepackage{amsmath, latexsym}
\usepackage[amsmath,thmmarks]{ntheorem}
\usepackage{algorithm}
\numberwithin{equation}{section}
\newtheorem{theorem}{Theorem}[section]
\newtheorem{lemma}{Lemma}[section]

\newtheorem{proposition}{Proposition}[section]
\newtheorem{remark}{Remark}[section]
\newtheorem{definition}{Definition}[section]
\newcommand{\beq}{\begin{equation}}
\newcommand{\eeq}{\end{equation}}
\newcommand{\beqn}{\begin{eqnarray}}
\newcommand{\eeqn}{\end{eqnarray}}
\date{}

\title{Max-min iterative methods and global error
bounds for the extended horizontal linear complementarity problem}
\author{Cui-Xia Li\thanks{lixiatkynu@126.com}, Shi-Liang Wu\thanks{Corresponding author: slwuynnu@126.com}, \\
{\small{\it$^{\ast, \dag}$School of Mathematics, Yunnan Normal University,}}\\
{\small{\it Kunming, Yunnan, 650500, P.R. China}}\\
{\small{\it$^{\ast, \dag}$Yunnan Key Laboratory of Modern Analytical Mathematics and Applications, }}\\
{\small{\it  Kunming, Yunnan, 650500, P.R. China}} }

\begin{document}

\maketitle

\begin{abstract}
To the best of our knowledge, since the extended horizontal linear
complementarity problem (EHLCP) was first introduced and studied by
Kaneko in 1977, no iterative methods or error analysis have been
developed for it due to the interdependence of its multiple unknowns
in a `chain-like' structure. This paper aims to address these gaps
by:
\begin{itemize}
\item proposing an equivalent fixed-point formulation of the EHLCP by
using a variable transformation technique with the max-min function;

\item developing the efficient max-min iterative methods for solving the EHLCP based on this fixed-point form,
along with their convergence analysis;

\item deriving global error bounds and computable estimates for the EHLCP.

\end{itemize}
Several numerical examples from applications such as multicommodity
market equilibrium and bilateral obstacle problems are given to
demonstrate the effectiveness of the proposed methods and bounds.
\end{abstract}



\pagestyle{myheadings} \thispagestyle{plain}
\markboth{\MakeUppercase{Shi-Liang Wu and Cui-Xia
Li}}{\MakeUppercase{MMIM and Global Error Bounds for the EHLCP}}

\section{Introduction}

For a given block matrix
\[
\mathrm{H}=(\mathrm{M},\mathrm{H}_{1},\ldots, \mathrm{H}_{m}) \
\mbox{with} \ \mathrm{M}, \mathrm{H}_{i}\in \mathbb{R}^{n\times n}
\]
and a given block vector
\[
\mathrm{d}=(\mathrm{q},\mathrm{d}_{1},\ldots, \mathrm{d}_{m-1})\
\mbox{with} \ \mathrm{q}\in \mathbb{R}^{n}\ \mbox{and}\
0<\mathrm{d}_{i}\in \mathbb{R}^{n} \ \mbox{for} \ i=1,2,\ldots,m-1,
\]
the problem of  finding $\mathrm{w},\mathrm{x}_{1}$, \ldots,
$\mathrm{x}_{m}\in \mathbb{R}^{n}$ such that
\addtocounter{equation}{1} \label{p11}
\begin{align}
&\mathrm{M}\mathrm{w}=\mathrm{q}+\sum_{i=1}^{m} \mathrm{H}_{i}\mathrm{x}_{i}, \label{1.1a} \tag{1.1a}\\
& \mathrm{w}, \mathrm{x}_{i}\geq0,i=1,\ldots,m, \label{1.1b} \tag{1.1b}\\
&\mathrm{w}^{T}\mathrm{x}_{1}=0,\label{1.1c} \tag{1.1c}\\
& \mathrm{x}_{i}\leq \mathrm{d}_{i},
(\mathrm{d}_{i}-\mathrm{x}_{i})^{T}\mathrm{x}_{i+1}=0,
i=1,2,\ldots,m-1, \label{1.1d} \tag{1.1d}
\end{align}
is said to be an extended horizontal linear complementarity problem
(EHLCP($\mathrm{H},\mathrm{d}$)), e.g., see \cite{Sznajder95} for
more details.  The EHLCP($\mathrm{H},\mathrm{d}$) covers some
classical complementarity problems, e.g.,

\begin{itemize}
\item taking $m=1$ in (1.1), the EHLCP($\mathrm{H},\mathrm{d}$) reduces
to the horizontal linear complementarity problem (HLCP($\mathrm{M},
\mathrm{H}_{1}, \mathrm{q}$)), see \cite{Fujisawa,Fujisawa2,Sun2};
\item taking $m=1$  and $\mathrm{H}_{1}=\mathrm{I}$ in (1.1), where,
hereafter, $\mathrm{I}$ denotes the identity matrix, the
EHLCP($\mathrm{H},\mathrm{d}$) reduces to the classical linear
complementarity problem (LCP($\mathrm{M}, \mathrm{q}$)), see
\cite{Cottle,Murty}.
\end{itemize}
It is known that many practical problems in scientific computing and
engineering applications, e.g., certain reinforced concrete frame
problem, a strictly convex quadratic program with bounded variables,
inventory theory and statistics, may reduce to the
EHLCP($\mathrm{H},\mathrm{d}$), see
\cite{Kaneko77,Kaneko80,Pang79,Kaneko79} and references therein.

For the EHLCP($\mathrm{H},\mathrm{d}$), Sznajder and Gowda in
\cite{Sznajder95} showed that the EHLCP($\mathrm{H},\mathrm{d}$) for
any $\mathrm{d}$ has a unique solution if and only if the block
matrix $\mathrm{H}$ has the column $\mathrm{W}$-property (which is
also
called as the $\mathrm{P}$-property in \cite{Kaneko78}). 
For $m=2$, from the point of view of least-element, Pang in
\cite{Pang79} established two least-element characterizations of the
solution of the EHLCP($\mathrm{H},\mathrm{d}$) with
$\mathrm{M}=\mathrm{H}_{2}=\mathrm{I}$. For the special matrix case,
Mezzadri and Galligani in \cite{Mezzadri4} generalized and analyzed
the concepts of diagonal dominance and irreducibility of column
representative matrices of a set, including the definitions of
particular sets of $\mathrm{M}$- and $\mathrm{H}$-matrices. In
addition, the finiteness, boundedness, nonemptyness, uniqueness and
connectedness issues of the solution set of the
EHLCP($\mathrm{H},\mathrm{d}$) have been investigated in recent
papers, see \cite{Yadav23,Yadav24}.

One important research problem for the complementarity problem (CP)
is the error bound because it often plays an important role in
theoretical analysis, including  sensitive analysis and verification
of the solutions, see \cite{Cottle}. For instance, assume that
$\mathrm{x}^{\ast}$ is the unique solution of the
LCP($\mathrm{M},\mathrm{q}$), a well-known global error bound for a
$\mathrm{P}$-matrix $\mathrm{M}$ was given by Mathias and Pang in
\cite{Mathias90}:
\begin{equation}\label{eq:12}
\|\mathrm{x}-\mathrm{x}^{\ast}\|_{\infty}\leq\frac{1+\|\mathrm{M}\|_{\infty}}{\alpha(\mathrm{M})}\|\mathrm{r}(\mathrm{x})\|_{\infty},
\ \forall \  \mathrm{x}\in\mathbb{R}^{n},
\end{equation}
where
\[
\mathrm{r}(\mathrm{x})=\min\{\mathrm{x},\mathrm{M}\mathrm{x}-\mathrm{q}\},
\alpha(\mathrm{M})=\min_{\|\mathrm{x}\|_{\infty}=1}\big\{\max_{1\leq
i\leq n}\mathrm{x}_{i}(\mathrm{M}\mathrm{x})_{i}\big\}.
\]
By using the equivalent form of the minimum function,  Chen and
Xiang in \cite{Chen06} obtained the following error bound in the
$p$-norm ($p\geq1, \mbox{or} \ p=\infty)$
\begin{equation}\label{eq:13}
\|\mathrm{x}-\mathrm{x}^{\ast}\|_{p}\leq\max_{d\in
[0,1]^{n}}\|(\mathrm{I}-\mathrm{D}+\mathrm{DM})^{-1}\|_{p}\|\mathrm{r}(\mathrm{x})\|_{p},
\forall \ \mathrm{x}\in \mathbb{R}^{n},
\end{equation}
where $\mathrm{D}=\mbox{diag}(d)$ with $d\in [0,1]^{n}$, which is
proved to be sharper than (\ref{eq:12}) for $p=\infty$. Very
recently, Wu and Wang \cite{Wu25} have extended the error bound
(\ref{eq:13}) to the extended vertical linear complementarity
problem (EVLCP). However, until now, the global error analysis of
EHLCP($\mathrm{H},\mathrm{d}$) has not been discussed, which
motivates us to investigate global error bounds in this paper.

As we know, both the linear complementarity problem and the extended
vertical linear complementarity problem essentially  involve only a
single unknown vector $\mathrm{x}$. Then with the aid of the
equivalent form of the minimum function, the global error bounds can
be deduced, see \cite{Chen06,Wu25}. However, the
EHLCP($\mathrm{H},\mathrm{d}$) contains multiple unknowns, which
depend on each other like a `chain', see (\ref{1.1c})-(\ref{1.1d}).
In such case, the technique used in \cite{Chen06,Wu25} is no longer
suitable for  establishing error bounds of the
EHLCP($\mathrm{H},\mathrm{d}$) because the constraint conditions
(\ref{1.1c})-(\ref{1.1d}) are different from the formula given by
the minimum function. So, the global error analysis of
EHLCP($\mathrm{H},\mathrm{d}$) has not been discussed so far.
Therefore, it needs to exploit the new technique to overcome these
disadvantages.

While the column $\mathrm{W}$-property of the system matrix
guarantees a unique solution to the EHLCP($\mathrm{H},\mathrm{d}$),
this theoretical result is of limited practical value without an
efficient numerical method to compute the solution. This motivates
our second key research objective: developing computationally
effective algorithms for solving the EHLCP($\mathrm{H},\mathrm{d}$).
For some special cases, there exist several methods. For instance,
when $m=2$ and $\mathrm{M}=\mathrm{H}_{2}=\mathrm{I}$ in (1.1), Pang
\cite{Pang79} proposed the modified Chandrasekaran's method for
$\mathrm{H}_{1}$
being a $\mathrm{Z}$-matrix.  
Kaneko and Pang \cite{Kaneko80} investigated  a parametric principal
pivoting algorithm to solve the LCP converted from (1.1) with
$\mathrm{M}=\mathrm{I}$. Both are direct methods. 
However, when employing these algorithms to solve the sparse and
large problem, it requires too many pivots, and can not keep the
sparsity structure of the complementarity system intact. Based on
this consideration, by using the relationship between  the linear
variational inequality problem (VIP) and the projection technique,
Ahn
\cite{Ahn} proposed a projection-type iteration method 
of $\mathrm{H}_{1}$ being an $\mathrm{H}_{+}$-matrix for $m=2$ and
$\mathrm{M}=\mathrm{H}_{2}=\mathrm{I}$ in (1.1) 
(see Method \ref{m33} in Section 3), which covers some existing
methods, e.g., block SOR method \cite{Cottle78}.

The development of numerical algorithms for solving the general
EHLCP($\mathrm{H},\mathrm{d}$) problems remains largely unexplored,
primarily due to the computational challenges posed by their chained
constraint structure. This paper addresses this gap by proposing
efficient iterative methods based on max-min function
transformations. Specifically, for the important special case where
$m=2$ with $\mathrm{M}=\mathrm{H}_{2}=\mathrm{I}$, we develop a
simple yet effective iterative algorithm when $\mathrm{H}_{1}$ is a
$\mathrm{P}$-matrix. Unlike direct methods that rely on pivoting
operations, our approach preserves the original problem structure
throughout iterations while maintaining sparsity. These features
make our methods particularly suitable for solving large-scale
sparse EHLCP($\mathrm{H},\mathrm{d}$) problems.

\vspace{0.3cm}

The key contributions of this work are:
\begin{itemize}
\item {\bf Equivalent fixed-point formulation:} A unified fixed-point
representation of the EHLCP is established by using a max-min
function and variable transformations, simplifying the problem to a
piecewise linear system.

\item {\bf Max-min iterative method:} A general max-min iterative method  (Method 3.1)
is proposed for the EHLCP, preserving matrix sparsity and avoiding
pivoting; a specialized method (Method 3.2) is developed for the
case where $m=2$ and $\mathrm{M}=\mathrm{H}_{2}=\mathrm{I}$  in
(1.1) and  $\mathrm{H}_{1}$ is a P-matrix, with proven convergence
under practical conditions.


\item {\bf Global error  bound:} The global error bounds of the EHLCP($\mathrm{H},\mathrm{d}$)
are derived, generalizing existing results for LCPs and HLCPs.
Computable bounds are provided for matrices with positive diagonal
dominance or strict diagonal dominance.

\item {\bf Numerical validation:} Experiments on real-world problems confirm the efficiency of the proposed methods and the sharpness of the error bounds.

\end{itemize}

The rest of this paper is organized as follows. In Section 2, the
equivalent form of the EHLCP($\mathrm{H},\mathrm{d}$) is exploited
by making use of variable transformation technique. In Section 3,
the max-min iterative method is discussed by using the above
equivalent fixed-point form, and a simple and efficient max-min
iterative method for $m=2$ and
$\mathrm{M}=\mathrm{H}_{2}=\mathrm{I}$ in (1.1)
is proposed. 
In Section 4, the useful framework of the global error bounds is
derived, and some computable estimates for the relevant upper bounds
are obtained under certain conditions. In Section 5, the numerical
results of some examples arising from different applications are
reported. Finally, in Section 6, we give some remarks to end this
paper.

The following notations,  definitions and lemmas \cite{Sznajder95,
Berman} will be used in the sequel.

Let $\mathbb{R}^{n\times n}$, $\mathbb{R}^{n}$ and $\mathbb{R}$ in
order denote the set of all $n\times n$ dimensional real matrices,
the set of all real $n$ dimensional vectors and the real field. We
use  small italic letters $x, y$, \dots, to represent real numbers,
lowercase letters $\mathrm{x}, \mathrm{y}$, \dots, to represent
vectors, capital letters $\mathrm{A}$, $\mathrm{B}$, \dots, to
represent matrices, and use calligraphic letters $\mathcal{A}$,
$\mathcal{B}$, \dots, to represent sets. For arbitrary positive
integer $n$, we denote $\mathcal{N}=\left \{1,\dots,n\right \}$.

For $\mathrm{x}, \mathrm{y}\in \mathbb{R}^{n}$,
$\min\left\{\mathrm{x}, \mathrm{y}\right\}$ represents the vector,
whose $i$th element is equal to \textrm{min}$\left\{x_{i},
y_{i}\right\}$; 
$\max\left\{\mathrm{x}, \mathrm{y}\right\}$ represents the vector,
whose $i$th element is equal to $\max\left\{x_{i}, y_{i}\right\}$.

Let $\mathrm{A} = (a_{ij})$ and $|\mathrm{A}|=(|a_{ij}|)$. A matrix
$\mathrm{A} =(a_{ij})$ is called a $\mathrm{Z}$-matrix if
$a_{ij}\leq0$ ($i\neq j$) for $i,j \in \mathcal{N}$; an
$\mathrm{M}$-matrix if $\mathrm{A}^{-1}\geq0$ and $a_{ij}\leq0$
($i\neq j$) for $i,j \in \mathcal{N}$; an $\mathrm{H}$-matrix if its
comparison matrix $\langle \mathrm{A}\rangle$ (i.e., $\langle
a\rangle_{ii}=|a_{ii}|, \langle a\rangle_{ij}=-|a_{ij}|$ $i\neq j$
for $i,j \in \mathcal{N}$) is an $\mathrm{M}$-matrix; an
$\mathrm{H}_{+}$-matrix if $\mathrm{A}$ is an $\mathrm{H}$-matrix
with $a_{ii} > 0$ for $i\in \mathcal{N}$; a column strictly diagonal
dominant (sdd) matrix if $|a_{ii}|>\Sigma_{j\neq i}|a_{ji}|,\ i\in
\mathcal{N}$; a row sdd
matrix if $|a_{ii}|>\Sigma_{j\neq i}|a_{ij}|,\ i\in \mathcal{N}$. 

Let $\mathrm{e}=(1,1,\ldots,1)^{T}$, $\rho(\cdot)$ be the
spectral radius of a matrix, and let $\mbox{diag}(d_{i})=\mbox{diag}(d_{1},...,d_n)$ and 
\begin{align*}
\mathcal{D}=\bigg\{(\mathrm{D}_{0},\mathrm{D}_{1},\ldots,\mathrm{D}_{m})~&|~\mathrm{D}_{i}=\mbox{diag}(d_{i1},d_{i2},\ldots,d_{in}), i=0,1,\ldots,m,\\
&\ d_{ij}\in [0,1]\ (j=1,\ldots,n)\ \mbox{and}\
\sum_{i=0}^{m}\mathrm{D}_{i}=\mathrm{I}\bigg\}.
\end{align*}
Without further illustration, the norm $\|\cdot\|$ for vectors means
the $p$-norm with $p \geq 1$, and for matrices it denotes the
corresponding induced operator norm. When other norms (e.g., the
spectral norm $\|\cdot\|_{2}$, the $1$-norm $\|\cdot\|_{1}$, or the
$\infty$-norm $\|\cdot\|_{\infty}$) are used, they will be
explicitly subscripted.

\begin{definition} \emph{\cite{Sznajder95}}
Let
$\mathcal{H}=\{\mathrm{M},\mathrm{H}_{1},\ldots,\mathrm{H}_{m}\}$
denote the set of matrices with
$\mathrm{M},\mathrm{H}_{1},\ldots,\mathrm{H}_{m}\in
\mathbb{R}^{n\times n}$. A matrix $\mathrm{R}\in \mathbb{R}^{n\times
n}$ is called a column representative of $\mathcal{H}$ if
\[
\mathrm{R}_{\cdot j}\in \{(\mathrm{M})_{\cdot
j},(\mathrm{H}_{1})_{\cdot j},\ldots, (\mathrm{H}_{m})_{\cdot j}\},
j=1,2,\ldots,n,
\]
where $\mathrm{S}_{\cdot j}$ denotes the j-th column of
$\mathrm{S}$.
\end{definition}

\begin{definition} \emph{\cite{Sznajder95}}
The set $\mathcal{H}$ is said to have the column
$\mathrm{W}$-property if  the determinants of all column
representative matrices of $\mathcal{H}$ are positive or negative.
\end{definition}

\begin{lemma} \label{y13} \emph{\cite{Sznajder95}} Let  $\mathcal{H}$ be given in Definition 1.1. The following statements are equivalent:
\begin{description}
\item $(i)$ $\mathcal{H}$ has the column $\mathrm{W}$-property;
\item $(ii)$ $ \mathrm{MX}_{0}+\mathrm{H}_{1}\mathrm{X}_{1}+\ldots+\mathrm{H}_{m}\mathrm{X}_{m}$
is nonsingular for arbitrary nonnegative diagonal matrices
$\mathrm{X}_{0}, \mathrm{X}_{1},\ldots, \mathrm{X}_{m}$ with
$\mbox{diag}(\mathrm{X}_{0}+\mathrm{X}_{1}+\ldots+\mathrm{X}_{m})>0$;
\item $(iii)$ the EHLCP$(\mathrm{H},\mathrm{d})$ has a unique solution for any vector $\mathrm{d}$.
\end{description}
\end{lemma}


\section{The equivalent form}
In this section, we will make use of variable transformation
technique to get an equivalent fixed-point form of the
EHLCP($\mathrm{H},\mathrm{d}$), which only involves one unknown.

Taking $m=1$ in (1.1), and
\begin{equation} \label{2.1}
\mathrm{w}=\max\left \{ \mathrm{0},-\mathrm{y}  \right \},
\mathrm{x}_{1}=\max\left \{ \mathrm{0},\mathrm{y} \right \} \
\mbox{with}\ \mathrm{y}\in \mathbb{R}^{n},
\end{equation}
it is easy to obtain an equivalent fixed-point form of
HLCP($\mathrm{M}, \mathrm{H}_{1}, \mathrm{q}$), also see \cite
{Bai10,Mezzadri2}. So, in the following we only need  to consider
the case of $m\geq2$ in (1.1). The following lemma is important to
design the iterative methods of EHLCP($\mathrm{H},\mathrm{d}$) and
get the global error bounds of EHLCP($\mathrm{H},\mathrm{d}$).

\begin{lemma}\label{y21}
Let $d_{0}=0$, $0<d_{i}\in \mathbb{R}$, $i=1,2,\ldots, m-1$ with
$m\geq2$. Then
\begin{align}
& w, x_{i}\geq0\ (i=1,\ldots,m), \label{2.2a} \tag{2.2a}\\
&wx_{1}=0,\label{2.2b} \tag{2.2b}\\
& x_{i}\leq d_{i}\ (i=1,2,\ldots,m-1), \label{2.2c} \tag{2.2c} \\
& (d_{i}-x_{i})x_{i+1}=0\  (i=1,2,\ldots,m-1). \label{2.2d}
\tag{2.2d}
\end{align}
if and only if there exists $y\in \mathbb{R}$ such that
\begin{equation*}
\begin{split}
&w=\max\left \{ 0,-y  \right \},\\
& x_{i}=\max\left \{ 0,\min\left \{ y-\sum_{j=0}^{i-1} d_{j},d_{i}
\right \} \right \}
\quad (i=1,2,\dots,m-1),\\
& x_{m}=\max\left \{ 0,y-\sum_{i=0}^{m-1} d_{i} \right \}.
\end{split}
\end{equation*}
Moreover, this result naturally extends, in component-wise form, to
vectors in $\mathbb{R}^{n}$.
\end{lemma}
\textbf{Proof.} $(\Leftarrow)$ Clearly, (\ref{2.2a})-(\ref{2.2c})
hold.

Now, we prove that  (\ref{2.2d}) holds. Firstly, we show 
%
 (\ref{2.2d})  for $i=1,2,\dots,m-2$, i.e.,
\begin{equation*}
\begin{split}
(d_{i}-x_{i})x_{i+1}=&\left(d_{i}-\max\left \{ 0,\min\left \{
y-\sum_{j=0} ^{i-1} d_{j},d_{i} \right \} \right
\}\right)\\
&  \times \max\left \{ 0,\min\left \{ y-\sum_{j=0}^{i} d_{j},d_{i+1}
\right \} \right \}=0.
\end{split}
\end{equation*}
Consider the following two cases:
\begin{description}
\item (1)
if $y-\sum_{j=0} ^{i-1} d_{j}\ge d_{i}$, then
\[
d_{i}-x_{i}=d_{i}-\max\left \{0,d_{i}\right \}=0,
\]
from which one may deduce $(x_{i}-d_{i})x_{i+1}=0$;
\item (2)
if $y-\sum_{j=0} ^{i-1} d_{j}< d_{i}$, i.e., $y-\sum_{j=0} ^{i}
d_{j}<0$, then
\begin{equation*}
\begin{split}
\max\left \{ 0,\min\left \{ y-\sum_{j=0}^{i}d_{j},d_{i+1} \right \}
\right \}=\max\left \{ 0, y-\sum_{j=0}^{i} d_{j} \right \}=0,
\end{split}
\end{equation*}
this implies that $(x_{i}-d_{i})x_{i+1}=0$.
\end{description}

Secondly, we prove that (\ref{2.2d}) holds  for $i=m-1$. By the
similar technique,  this case also holds, whose proof is omitted.

$(\Rightarrow)$ We consider two main cases on the base of the values
of  $x_i$.

Case 1: All  $x_i = 0$  for  $i = 1, \ldots, m$.

From (2.2b),  $w x_1 = 0$  and  $x_1 = 0$, so  $w$  can be any
non-negative number. We need to find  $y$  such that:
\begin{equation*}
\begin{split}
&w=\max\left \{ 0,-y  \right \},\\
& for \ i = 1, \ldots, m-1 , x_i=\max\left \{ 0,\min\left \{
y-\sum_{j=0}^{i-1} d_{j},d_{i} \right \} \right \}=0,\\
& x_{m}=\max\left \{ 0,y-\sum_{i=0}^{m-1} d_{i} \right \}=0.
\end{split}
\end{equation*}
Since  $d_0 = 0$  and  $d_i > 0$  for  $i \geq 1$, we choose  $y
\leq 0 $. Then for any  $i \geq 1$:
\begin{description}
\item $\bullet$ $y - \sum_{j=0}^{i-1} d_j \leq y \leq 0$, so  $\min\left\{ y - \sum_{j=0}^{i-1} d_j, d_i \right\} \leq
0$ ($d_i > 0$), hence \[\max\left\{ 0,  \min\left\{y -
\sum_{j=0}^{i-1} d_j, d_i  \right\} \right\} = 0.\]
\item $\bullet$ Similarly,  $y - \sum_{i=0}^{m-1} d_i \leq 0$, so  $\max\left\{ 0, y
- \sum_{i=0}^{m-1} d_i \right\} = 0$.
\end{description}
Of course, $w = \max\{0, -y\}\geq 0$ for $y \leq 0 $. Thus, all
expressions are satisfied.
%
%

\vspace{0.3cm}

Case 2: Some  $x_i > 0.$

Let \(k\) be the largest index such that \(x_k > 0\) (so \(x_i = 0\)
for all \(i > k\)).

We first show that \(x_1 > 0\). If \(x_1 = 0\), then from (2.2d)
with \(i = 1\), we have \((d_1 - x_1)x_2 = d_1 x_2 = 0\). Since
\(d_1 > 0\), this implies \(x_2 = 0\). Proceeding inductively, we
obtain \(x_i = 0\) for all \(i \ge 1\), which contradicts the
assumption that some \(x_i > 0\). Hence, \(x_1 > 0\). Further,  from
(2.2b), we have $w = 0$.



Now we consider two subcases based on the value of \(k\):

%
\begin{description}
\item
$\bullet$ \(k = m\) (i.e., \(x_m > 0\)).
          Then from (2.2d) for \(i = m-1\), we get \((d_{m-1} - x_{m-1})x_m = 0\), which forces \(x_{m-1} = d_{m-1}\) since \(x_m > 0\).
Next, from (2.2d) with \(i = m-2\), we have \((d_{m-2} -
x_{m-2})x_{m-1} = 0\). Since \(x_{m-1} = d_{m-1} > 0\), this forces
\(x_{m-2} = d_{m-2}\).
          Repeating this backward induction yields \(x_i = d_i\) for all \(i = 1, \ldots, m-1\).


%

\item
$\bullet$ If  $k < m$, then  $x_i = 0$  for  $i > k$.
          We show by backward induction that \(x_i = d_i\) for all \(i = 1, \ldots, k-1\).
          \begin{itemize}
              \item For \(i = k-1\): From (2.2d), \((d_{k-1} - x_{k-1})x_k = 0\). Since \(x_k > 0\), we must have \(d_{k-1} - x_{k-1} = 0\), i.e., \(x_{k-1} = d_{k-1}\).
              \item Now assume \(x_j = d_j\) for some \(j\) with \(1 < j \le k-1\). From (2.2d) with \(i = j-1\), we have \((d_{j-1} - x_{j-1})x_j = 0\). Since \(x_j = d_j > 0\), this forces \(x_{j-1} = d_{j-1}\).
          \end{itemize}
          By induction, we obtain \(x_i = d_i\) for all \(i = 1, \ldots, k-1\). Since \(k < m\) implies \(k \le m-1\), condition (2.2c) applies to \(x_k\), giving \(x_k \le d_k\).


\end{description}

Now, construct  $y$:

\begin{description}
\item Subcase 2a:  $k = m$ (i.e.,  $x_m > 0$). Let $y =
\sum_{j=1}^{m-1} d_j + x_m$.
\item  Subcase 2b:  $k <m$. Let  $y =
\sum_{j=1}^{k-1} d_j + x_k$  (if  $k = 1$, then  $\sum_{j=1}^{0} d_j
= 0$, so \( y = x_1 \)).
\end{description}

We now verify the conclusions:

\begin{description}
\item Since  $w = 0$  and  $y \geq 0$, we have $\max\{0, -y\} = 0 = w$.

\item  For $x_i$ ($i = 1, \ldots, m-1$):

\begin{description}
\item Subcase 2a:  $y = \sum_{j=1}^{m-1} d_j + x_m$.

       For $i = 1, \ldots, m-1$,

\[y - \sum_{j=0}^{i-1} d_j = \sum_{j=1}^{m-1} d_j + x_m -
\sum_{j=1}^{i-1} d_j = \sum_{j=i}^{m-1} d_j + x_m \geq d_i\ \
(\sum_{j=i}^{m-1} d_j \geq d_i).
\]  Thus, $\min\left\{ y - \sum_{j=0}^{i-1} d_j, d_i
\right\} = d_i$, and $\max\left\{0, d_i \right\} = d_i = x_i$. I.e.,
\[
 x_i=\max\left\{0, \min\left\{ y - \sum_{j=0}^{i-1} d_j, d_i
\right\} \right\}.
\]
\item Subcase 2b:  $y = \sum_{j=1}^{k-1} d_j + x_k$.
\begin{description}
\item If  $i < k$, then
 \[y - \sum_{j=0}^{i-1} d_j = \sum_{j=1}^{k-1} d_j + x_k - \sum_{j=1}^{i-1} d_j = \sum_{j=i}^{k-1} d_j + x_k \geq
 d_i\ \ (\sum_{j=i}^{k-1} d_j \geq d_i). \] So  $\min\left\{ y -
\sum_{j=0}^{i-1} d_j, d_i \right\} = d_i$,  and $\max\left\{0, d_i
\right\} = d_i = x_i$. Further,
\[
 x_i=\max\left\{0, \min\left\{ y -
\sum_{j=0}^{i-1} d_j, d_i \right\} \right\}.
\]
\item If $i = k$, then $y - \sum_{j=0}^{k-1} d_j = x_k$.
      So  $\min\left\{ x_k, d_k \right\} = x_k$ ($x_k \leq d_k$),
            and $\max\left\{0, x_k \right\} = x_k$. I.e.,
            \[
x_k=\max\left\{0, \min\left\{ y - \sum_{j=0}^{k-1} d_j, d_k \right\}
\right\}.
            \]

\item If  $i > k$, then
      \[y - \sum_{j=0}^{i-1} d_j = \sum_{j=1}^{k-1} d_j + x_k - \sum_{j=1}^{i-1}
      d_j.\]
      Since $i-1 \geq k$, we have  $\sum_{j=1}^{i-1} d_j \geq \sum_{j=1}^{k} d_j
      $. So,
     \[ y - \sum_{j=0}^{i-1} d_j \leq \sum_{j=1}^{k-1} d_j + x_k - \sum_{j=1}^{k} d_j = x_k - d_k \leq
     0.\]
      Thus, $\min\left\{ y - \sum_{j=0}^{i-1} d_j, d_i \right\} \leq 0$,
            and
            \[\max\left\{0, \min\left\{ y - \sum_{j=0}^{i-1} d_j, d_i \right\} \right\} = 0 =
            x_i.\]
\end{description}
\end{description}
For $x_m$:
\begin{description}
\item  Subcase 2a: $y - \sum_{i=0}^{m-1} d_i = \sum_{j=1}^{m-1} d_j + x_m - \sum_{j=0}^{m-1} d_j =
x_m$ ($\sum_{j=0}^{m-1} d_j = \sum_{j=1}^{m-1} d_j $ since $d_0 =
0$). So $\max\left\{0, x_m \right\} = x_m$, i.e.,
\[
x_m=\max\left\{0, y - \sum_{i=0}^{m-1} d_i \right\}.
\]
\item  Subcase 2b:
        Since $k < m$,  $x_m = 0$.
    \[y - \sum_{i=0}^{m-1} d_i = \sum_{j=1}^{k-1} d_j + x_k - \sum_{j=0}^{m-1} d_j=x_k-( d_k+\ldots+d_{m-1}) \leq x_k - d_k \leq
    0.\]
    So,  $\max\left\{0, y - \sum_{i=0}^{m-1} d_i \right\} = 0 = x_m$.
\end{description}
\end{description}
In all cases, the constructed  $y$  satisfies the required
expressions. $\hfill{} \Box$



\begin{remark} \label{r21} 
Lemma \ref{y21} is also true for $m=1$. In this case,  both (2.2a)
and (2.2b) remain valid. (\ref{2.2c}) makes no sense since no
$x_{0}$ exists.
\end{remark}

By Lemma \ref{y21} and Remark \ref{r21}, it is easy to prove that
the EHLCP($\mathrm{H},\mathrm{d}$) can be transformed to the
following unified and equivalent fixed-point formula. Concretely,
see Proposition \ref{pro21}.

\begin{proposition}\label{pro21} Let $\mathrm{y}\neq0\in \mathbb{R}^{n}$. For the EHLCP$(\mathrm{H},\mathrm{d})$, setting
\begin{equation} \label{2.3} \tag{2.3}
\begin{split}
&\mathrm{w}=\max\left \{ \mathrm{0},-\mathrm{y}  \right \},\\
& \mathrm{x}_{i}=\max\left \{ 0,\min\left \{
\mathrm{y}-\sum_{j=0}^{i-1} \mathrm{d}_{j},\mathrm{d}_{i} \right \}
\right \}
\quad (i=1, 2,\dots,m-1),\\
& \mathrm{x}_{m}=\max\left \{ 0,\mathrm{y}-\sum_{i=0}^{m-1}
\mathrm{d}_{i} \right \},
\end{split}
\end{equation}
where  $\mathrm{d}_{0}=0$. Then the EHLCP$(\mathrm{H},\mathrm{d})$
is equal to finding $\mathrm{y}$ such that
\begin{align}\label{2.4} \tag{2.4}
\mathrm{M}\max\{0,-\mathrm{y}\}=&\mathrm{q}+\sum_{i=1}^{m-1}
\mathrm{H}_{i}\max\left \{ 0,\min\left \{ \mathrm{y}
-\sum_{j=0}^{i-1} \mathrm{d}_{j},\mathrm{d}_{i} \right \} \right \}\\
&+\mathrm{H}_{m}\max\left \{ 0,\mathrm{y}-\sum_{i=0}^{m-1}
\mathrm{d}_{i} \right \}.\nonumber
\end{align}
\end{proposition}

\begin{remark}In Proposition \ref{pro21}, $\mathrm{y}\neq0\in \mathbb{R}^{n}$ is
reasonable.  The reason is that if $\mathrm{y}=0$, then Eq.
(\ref{2.4}) is not valid for $\mathrm{q}\neq0$. Such a case,
$\mathrm{w}=\mathrm{x}_{1}=\ldots=\mathrm{x}_{m}=0$ is not the
solution of the EHLCP$(\mathrm{H},\mathrm{d})$ for
$\mathrm{q}\neq0$. Hence, in the later discussion, without further
illustration, we always assume that $\mathrm{y}\neq0\in
\mathbb{R}^{n}$.
\end{remark}

In the sequel, Eq. (\ref{2.4}) will be employed to establish the
max-min iterative methods and the global error bounds for the
EHLCP($\mathrm{H},\mathrm{d}$).

\begin{remark} Lemma \ref{y21} still holds if we introduce a positive diagonal matrix
$\mathrm{\Omega}\in \mathbb{R}^{n\times n}$ for $\mathrm{w}$ and a
positive diagonal matrix $\mathrm{\Delta}\in \mathbb{R}^{n\times n}$
for $\mathrm{x}_{m}$, respectively, i.e.,
\[
\mathrm{w}=\mathrm{\Omega}\max\{ 0,-\mathrm{y} \} \mbox{ and}\
\mathrm{x}_{m}=\mathrm{\Delta} \max\left \{
0,\mathrm{y}-\sum_{i=0}^{m-1} \mathrm{d}_{i} \right \}.
\]
\end{remark}

Combing Lemma \ref{y21} with Remark 2.3, the equivalent fixed-point
formula of EHLCP($\mathrm{H},\mathrm{d}$) with parameters can also
be given, which is omitted.

\section{The max-min iterative method}
In this section,  our main task is to propose the max-min iterative
methods. Clearly, Eq. (\ref{2.4}) is a piecewise linear systems
(PLS) that only involves one unknown vector $\mathrm{y}$.
If Eq. (\ref{2.4}) has a unique solution
$\mathrm{y}$, then the solution $\mathrm{w},
\mathrm{x}_{1},\ldots, \mathrm{x}_{m}$ of the EHLCP($\mathrm{H,d}$)
can be directly obtained from by  (\ref{2.3}).

By Lemma \ref{y13}, the column $\mathrm{W}$-property of $\mathcal{H}$ guarantees that the EHLCP$(\mathrm{H}, \mathrm{d})$ has a unique solution for any vector $\mathrm{d}$.
Moreover, since the transformation (\ref{2.3}) establishes an equivalence between the EHLCP and the PLS (\ref{2.4}), it follows that under the same column $\mathrm{W}$-property, Eq. (2.4) also admits a unique solution $\mathrm{y}^*$.
This equivalence not only validates our variable transformation but also ensures that any iterative scheme designed for (\ref{2.4}) targets a well-posed problem.

To facilitate the design of iterative methods, we first establish some useful properties of the max-min function.
\begin{lemma}\label{y31}
Let $a,b\in \mathbb{R}$ with $a>0$. Then for any $x,y\in \mathbb{R}$
with $x\neq y$, there exists $\alpha\in [0,1]$ such that
\begin{equation*}
\max\{0,y-a\}-\max\{0,x-a\}=\alpha(y-x).
\end{equation*}
Moreover, if $0<a<b$, then there exist $\alpha$ and $\beta\in [0,1]$
such that  $0\leq \beta\leq \alpha\leq1$ and
\begin{equation*}
(\max\{0,y-a\}-\max\{0,x-a\})-(\max\{0,y-b\}-\max\{0,x-b\})=(\alpha-\beta)(y-x).
\end{equation*}
Furthermore, these results are also true for component-wise form to
vectors in $\mathbb{R}^{n}$.
\end{lemma}

\textbf{Proof.} Its proof is straightforward, which is omitted.

\vspace{0.3cm}

By Lemma \ref{y31}, 
we have

\begin{lemma} \label{y32} Let  $h_{i}\in \mathbb{R}\ (i=1,\ldots,m)$,  $d_{i}\in \mathbb{R}\ (i=0,1,\ldots,m-1)$ be defined as in Lemma \ref{y21}. Then
there are $\lambda_{i}\geq0,\sum_{i=0}^{m}\lambda_{i}=1$ such that
for $y,y^{\ast}\in \mathbb{R}$,
\begin{align*}
&\sum_{i=1}^{m-1} h_{i}\max \bigg\{ 0,\min \bigg \{
y-\sum_{j=0}^{i-1} d_{j},d_{i}  \bigg\}  \bigg\}
+h_{m}\max \bigg \{ 0,y-\sum_{i=0}^{m-1} d_{i}\bigg\}-h_{0}\max\{0,-y\} \nonumber\\
&-\bigg(\sum_{i=1}^{m-1} h_{i}\max \bigg\{ 0,\min \bigg\{
y^{\ast}-\sum_{j=0}^{i-1} d_{j},d_{i} \bigg \} \bigg\}\\
&+h_{m}\max \bigg\{ 0,y^{\ast}-\sum_{i=0}^{m-1} d_{i} \bigg\}-h_{0}\max\{0,-y^{\ast}\} \bigg)\\
=&(h_{0}\lambda_{0}+h_{1}\lambda_{1}+\ldots+h_{m}\lambda_{m})(y-y^{\ast}).
\end{align*}
Moreover, this result naturally extends, in component-wise form, to
vectors in $\mathbb{R}^{n}$.
\end{lemma}
\textbf{Proof.} By making use of
\begin{align}\label{3.1}
\max\{0,\min\{c, x\}\}=\max\{0,x\}-\max\{0,x-c\}\  \mbox{for} \ c>0
\end{align}
and
\begin{equation}\label{3.2}
\max\{0,-\mathrm{y}\}+\mathrm{y}=\max\{0,\mathrm{y}\},
\end{equation}
we have
\begin{align*}
&\sum_{i=1}^{m-1} h_{i}\max \bigg\{ 0,\min \bigg\{ y-\sum_{j=0}^{i-1} d_{j},d_{i} \bigg\}  \bigg\}+h_{m}\max \bigg\{ 0,y-\sum_{i=0}^{m-1} d_{i}\bigg\}-h_{0}\max\{0,-y\}\\
&-\bigg(\sum_{i=1}^{m-1} h_{i}\max \bigg\{ 0,\min \bigg\{
y^{\ast}-\sum_{j=0}^{i-1} d_{j},d_{i}\bigg\} \bigg\} +h_{m}\max
\bigg\{ 0,y^{\ast}-\sum_{i=0}^{m-1} d_{i}
\bigg\}-h_{0}\max\{0,-y^{\ast}\} \bigg)\\
=&\sum_{i=1}^{m} h_{i}\max\bigg\{0,y-\sum_{i=1}^{m} d_{i-1}\bigg\}-
\sum_{i=1}^{m-1} h_{i}\max\bigg\{0,y-\sum_{i=1}^{m} d_{i}\bigg\}+h_{0}y-h_{0}\max\{0,y\}\\
\end{align*}
\begin{align*}
&-\bigg(\sum_{i=1}^{m} h_{i}\max\bigg\{0,y^{\ast}-\sum_{i=1}^{m}
d_{i-1}\bigg\}-
\sum_{i=1}^{m-1} h_{i}\max\bigg\{0,y^{\ast}-\sum_{i=1}^{m} d_{i}\bigg\}+h_{0}y^{\ast}-h_{0}\max\{0,y^{\ast}\}\bigg)\\
=&h_{0}(y-y^{\ast})+\sum_{i=1}^{m}(h_{i}-h_{i-1})\bigg(\max\bigg\{0,y-\sum_{i=1}^{m}d_{i-1}\bigg\}-\max\bigg\{0,y^{\ast}-\sum_{i=1}^{m}d_{i-1}\bigg\}\bigg)\\
=&h_{0}(y-y^{\ast})+\sum_{i=1}^{m}(h_{i}-h_{i-1})\nu_{i}(y-y^{\ast})\\
=&(h_{0}(1-\nu_{1})+h_{1}(\nu_{1}-\nu_{2})+h_{2}(\nu_{2}-\nu_{3})+\ldots+h_{m-1}(\nu_{m-1}-\nu_{m})+h_{m}\nu_{m})(y-y^{\ast})\\
=&(h_{0}\lambda_{0}+h_{1}\lambda_{1}+h_{2}\lambda_{2}+h_{3}\lambda_{3}+\ldots+h_{m-1}\lambda_{m-1}+h_{m}\lambda_{m})(y-y^{\ast}),
\end{align*}
where the third to the last equalities follow from Lemma \ref{y31},
i.e., there exists $0\leq\nu_{i}\leq1 \ (i=1,2,\ldots,m)$ such that
it holds, and
\[
\lambda_{0}=1-\nu_{1},\lambda_{1}=\nu_{1}-\nu_{2},\lambda_{3}=\nu_{2}-\nu_{3},\ldots,
\lambda_{m-1}=\nu_{m-1}-\nu_{m},\lambda_{m}=\nu_{m}.
\]
From Lemma \ref{y31}, it is easy to check that
\[
1\geq \nu_{1}\geq \nu_{2}\geq \nu_{3}\geq \cdots \geq \nu_{m}\geq0,
\]
and
\[
\lambda_{0}+\lambda_{1}+\ldots+\lambda_{m}=1.
\]
This completes proof. $\hfill{} \Box$

As a direct consequence of Proposition \ref{pro21}, we obtain the following fundamental identity that reveals the intrinsic relationship among the variables.

\begin{theorem}\label{th32}
Let $\mathrm{y}\in \mathbb{R}^{n}$, $\mathrm{w}(\mathrm{y}),
\mathrm{x}_{1}(\mathrm{y}),\ldots, \mathrm{x}_{m}(\mathrm{y})$ be
given by (\ref{2.3}), i.e.,
\begin{align*}
\mathrm{w}(\mathrm{y}) &= \max\{0,-\mathrm{y}\}, \\
\mathrm{x}_{i}(\mathrm{y}) &= \max\Bigg\{0,\min\Bigg\{\mathrm{y}-\sum_{j=0}^{i-1}\mathrm{d}_{j},\mathrm{d}_{i}\Bigg\}\Bigg\} \quad (i=1,\dots,m-1), \\
\mathrm{x}_{m}(\mathrm{y}) &=
\max\left\{0,\mathrm{y}-\sum_{i=0}^{m-1}\mathrm{d}_{i}\right\}.
\end{align*} Then
\begin{equation}\label{n31}
\sum_{i=1}^{m}\mathrm{x}_{i}(\mathrm{y})=\mathrm{y}+\mathrm{w}(\mathrm{y}).
\end{equation}
\end{theorem}

\textbf{Proof.} Obviously. $\hfill{} \Box$ 

\vspace{0.3cm}

In view of Eq. (\ref{3.2}), now we may propose the following max-min
iterative iterative method:

\textbf{Method 3.1} \label{m31} Let $\mathcal{H}$ have the column
$\mathrm{W}$-property, and let $\mathrm{y}^{0}\in \mathbb{R}^{n}$ be
an initial guess.
For $k=0,1,2,\ldots$, until convergence, compute $\mathrm{y}^{k+1}$
by
\begin{align}\label{3.4}
\mathrm{My}^{k+1}=\mathrm{M}\max\{0,\mathrm{y}^{k}\}-\mathrm{q}-\phi(\mathrm{y}^{k}),
\end{align}
where
\begin{align*}
\phi(\mathrm{y}^{k})=\sum_{i=1}^{m-1} \mathrm{H}_{i}\max\Bigg\{
0,\min\Bigg \{ \mathrm{y}^{k}-\sum_{j=0}^{i-1}
\mathrm{d}_{j},\mathrm{d}_{i} \Bigg \} \Bigg \}
+\mathrm{H}_{m}\max\left \{ 0,\mathrm{y}^{k}-\sum_{i=0}^{m-1}
\mathrm{d}_{i} \right \}.
\end{align*}
Let $\mathrm{y}^{*}$ be the convergence point of $\mathrm{y}^{k}$.
Then set
\begin{align*}
&\mathrm{w}^{*}=\max\left \{ \mathrm{0},-\mathrm{y}^{*}  \right \},\\
& \mathrm{x}^{*}_{i}=\max\Bigg\{ 0,\min\Bigg\{
\mathrm{y}^{*}-\sum_{j=0}^{i-1} \mathrm{d}_{j},\mathrm{d}_{i}\Bigg
\} \Bigg\}
\quad (i=1, 2,\dots,m-1),\\
& \mathrm{x}^{*}_{m}=\max\left \{ 0,\mathrm{y}^{*}-\sum_{i=0}^{m-1}
\mathrm{d}_{i} \right \}.
\end{align*}

Clearly,  by Proposition \ref{pro21}, the vectors $\mathrm{w}^{*},
\mathrm{x}^{*}_i$ generated by Method 3.1 are the solution of
EHLCP($\mathrm{H,d}$).

It is also known that the proposed method, Method 3.1, does not
involve any pivoting operations as in the direct method.
A key advantage of Method 3.1 (and similarly, Method 3.2) is its
inherent ability to preserve sparsity. Each iteration involves only
two types of operations:
\begin{enumerate}
\item Matrix-vector multiplications with the original matrices $\mathrm{M}, \mathrm{H}_1, \ldots, \mathrm{H}_m$, and
\item Element-wise vector operations ($\max$, $\min$, addition, subtraction).
\end{enumerate}
For large-scale sparse matrices, the matrix-vector products are
computationally efficient and do not incur fill-in (unlike direct
methods that require matrix factorization or pivot operations). The
vectors $\mathrm{y}^k$ and the results of the element-wise
operations remain dense, but their storage requirements grow only
linearly with $n$. This structural preservation makes the method
particularly suitable for solving large-scale sparse
EHLCP($\mathrm{H,d}$) arising in practical applications such as
those presented in Section 5.

\vspace{0.3cm}

The following result is the global convergence theorem for Method
3.1.

\begin{theorem}\label{th33}  Let $\mathcal{H}$ have the column
$\mathrm{W}$-property.
If
\begin{equation}\label{3.5}
\rho\bigg(\sum_{i=1}^{m}|\mathrm{I}-\mathrm{M}^{-1}\mathrm{H}_{i}|\bigg)<1\
\mbox{or}\
\sum_{i=1}^{m}\|\mathrm{I}-\mathrm{M}^{-1}\mathrm{H}_{i}\|<1,
\end{equation}
then the iterative sequence $\{\mathrm{y}^{k}\}^{+\infty}
_{k=0}\subset\mathbb{R}^{n}$ generated by Method 3.1 converges to
the unique solution $\mathrm{y}^{\ast}\in \mathbb{R}^{n}$ of the PLS
(\ref{2.4}) for an arbitrary initial vector
$\mathrm{y}^{0}\in\mathbb{R}^{n}$.
\end{theorem}

\textbf{Proof.} From (\ref{2.4}), $\mathrm{y}^{\ast}$ satisfies
\begin{align} \label{3.6}
0=&-\Bigg(\sum_{i=1}^{m-1} \mathrm{H}_{i}\max\left \{ 0,\min\left \{
\mathrm{y}^{\ast}
-\sum_{j=0}^{i-1} \mathrm{d}_{j},\mathrm{d}_{i} \right \} \right \}\\
&+\mathrm{H}_{m}\max\left \{ 0,\mathrm{y}^{\ast}-\sum_{i=0}^{m-1}
\mathrm{d}_{i} \right
\}-\mathrm{M}\max\{0,-\mathrm{y}^{\ast}\}\Bigg)-\mathrm{q}.\nonumber
\end{align}
And Eq. (\ref{3.4}) can be expressed as
\begin{align}\label{3.7}
\mathrm{My}^{k+1}-\mathrm{My}^{k}=&-\Bigg(\sum_{i=1}^{m-1}
\mathrm{H}_{i}\max\left \{ 0,\min\left \{
\mathrm{y}^{k}-\sum_{j=0}^{i-1}
\mathrm{d}_{j},\mathrm{d}_{i} \right \} \right \}\\
&+\mathrm{H}_{m}\max\left \{ 0,\mathrm{y}^{k}-\sum_{i=0}^{m-1}
\mathrm{d}_{i} \right
\}-\mathrm{M}\max\{0,-\mathrm{y}^{k}\}\Bigg)-\mathrm{q}. \nonumber
\end{align}
Subtracting  (\ref{3.6}) from (\ref{3.7}), from Lemma \ref{y32},
there exist nonnegative diagonal matrices
$\mathrm{D}^{k}_{i}=\mbox{diag}(d^{k}_{i})$ with  $d^{k}_{i}\in [0,1]^{n}$
$(i=0,1,\ldots,m)$ and $\sum_{i=0}^{m}\mathrm{D}^{k}_{i}=\mathrm{I}$,
i.e., $(\mathrm{D}^{k}_{0},\mathrm{D}^{k}_{1},\ldots,\mathrm{D}^{k}_{m})\in
\mathcal{D}$ such that
\begin{align*}
\mathrm{My}^{k+1}-\mathrm{My}^{k}=
-(\mathrm{M}\mathrm{D}^{k}_{0}+\mathrm{H}_{1}\mathrm{D}^{k}_{1}+\ldots+\mathrm{H}_{m}\mathrm{D}^{k}_{m})(\mathrm{y}^{k}-\mathrm{y}^{\ast}).
\end{align*}

Since the block matrix $\mathrm{H}$ has the column
$\mathrm{W}$-property,  $\mathrm{M}$ is nonsingular. 
Hence we have
\begin{align*}
\mathrm{y}^{k+1}-\mathrm{y}^{\ast}&=(\mathrm{I}-\mathrm{M}^{-1}(\mathrm{M}\mathrm{D}^{k}_{0}+\mathrm{H}_{1}\mathrm{D}^{k}_{1}
+\ldots+\mathrm{H}_{m}\mathrm{D}^{k}_{m}))(\mathrm{y}^{k}-\mathrm{y}^{\ast}),\\
&=\bigg(\sum_{i=1}^{m}(\mathrm{I}-\mathrm{M}^{-1}\mathrm{H}_{i})D^{k}_{i}\bigg)(\mathrm{y}^{k}-\mathrm{y}^{\ast}),
\end{align*}
which together with the condition (\ref{3.5})  the iterative series
$\{\mathrm{y}^{k}\}^{+\infty} _{k=0}\subset\mathbb{R}^{n}$ given by
Method 3.1 converge. This proves Theorem \ref{th33}. $\hfill{} \Box$

\vspace{0.3cm}

It is noted that from PLS (\ref{2.4}) and Remark 2.3, some other
iterative methods such as the parameter version of Method 3.1, the
Newton-type method, the modulus-type method can be established,
which remains the discussion in the future.

\vspace{0.3cm}

Next, following the methodology of Method 3.1, we can design a
simple and efficient iterative  method for solving the following
special EHLCP with $m=2,~\mathrm{M}=\mathrm{H}_2=\mathrm{I}$, i.e.,
\begin{equation} \label{3.9}
\begin{split}
&\mathrm{w}=\mathrm{q}+ \mathrm{H}_{1}\mathrm{x}_{1}+\mathrm{x}_{2}, \\
& \mathrm{w}, \mathrm{x}_{1}, \mathrm{x}_{2}\geq0, \mathrm{b}>0, \\
&\mathrm{w}^{T}\mathrm{x}_{1}=0,\\
& \mathrm{x}_{1}\leq \mathrm{b},
(\mathrm{b}-\mathrm{x}_{1})^{T}\mathrm{x}_{2}=0,
\end{split}
\end{equation}
which is often from many practical problems, e.g., multicommodity
market equilibrium problem, quadratic programming (QP) problem, the
bilateral obstacle problem and others. Assuming that
$\mathrm{H}_{1}$ is a $\mathrm{P}$-matrix. Then by Lemma \ref{y13}
the EHLCP($\mathrm{H},\mathrm{d}$) (\ref{3.9}) has a unique
solution.

Based on Lemma \ref{y21} and Remark 2.2, we take
\[
\mathrm{w}=\Omega\max\{0,-\mathrm{y}\}, \mathrm{x}_{1}=\max\left \{
0,\min\left \{ \mathrm{y},\mathrm{b} \right \} \right \},
\mathrm{x}_{2}=\Delta\max\{0,\mathrm{y}-\mathrm{b}\}
\]
for the EHLCP($\mathrm{H},\mathrm{d}$) (\ref{3.9}), where $\Delta$
and $\mathrm{\Omega}$ are any positive diagonal matrices. Here,  for
the sake of simplicity, we only consider this case where $\Omega=
\Delta$. In such case, its PLS is of the form
\begin{equation*}
\Omega\max\{0,-\mathrm{y}\}=\mathrm{q}+\mathrm{H}_{1}\max\left \{
0,\min\left \{ \mathrm{y},\mathrm{b} \right \} \right
\}+\Omega\max\left \{0,\mathrm{y}-\mathrm{b} \right \},
\end{equation*}
which together with (\ref{3.1}) and (\ref{3.2}) gives
\begin{equation} \label{3.10}
\Omega\mathrm{y}=-((\mathrm{H}_{1}-\Omega)\max\left \{ 0,\min\left
\{ \mathrm{y},\mathrm{b} \right \} \right \}+\mathrm{q}).
\end{equation}

Based on  Eq. (\ref{3.10}), we can  design the following max-min
iterative method for solving the EHLCP($\mathrm{H},\mathrm{d}$)
(\ref{3.9}).


\textbf{Method 3.2 }\label{m32}  Given an arbitrary initial guess
$\mathrm{y}^{0}\in \mathbb{R}^{n}$, for $k=0,1,2,\ldots$, until
convergence, compute $\mathrm{y}^{k+1}$ by
\begin{equation}\label{3.11}
\Omega\mathrm{y}^{k+1}=-((\mathrm{H}_{1}-\Omega)\max\left \{
0,\min\left \{ \mathrm{y}^{k},\mathrm{b} \right \} \right
\}+\mathrm{q}), k=0, 1, \ldots,
\end{equation}
where $\mathrm{\Omega}$ is any positive diagonal matrix. Let
$\mathrm{y}^{*}$ be the convergence point of $\mathrm{y}^{k}$. Then
set
\begin{align*}
\mathrm{w}^{*}&=\Omega\max\{0,-\mathrm{y}^{*}\},\\
 \mathrm{x}_{1}^{*}&=\max\left
\{ 0,\min\left \{ \mathrm{y}^{*},\mathrm{b} \right \} \right \},\\
\mathrm{x}_{2}^{*}&=\Omega\max\{0,\mathrm{y}^{*}-\mathrm{b}\}.
\end{align*}

It is easy to see that the vectors $\mathrm{w}^{*},
\mathrm{x}^{*}_1$ and $\mathrm{x}^{*}_2$  generated by Method 3.2
are the numerical solution of EHLCP($\mathrm{H},\mathrm{d}$)
(\ref{3.9}).

Similarly, Method 3.2 does not involve any pivoting operations,
which implies that it is different from the modified
Chandrasekaran's method \cite{Pang79}. Clearly, in the proposed
algorithm, the original data of matrices and vectors always remain
unchanged in each iteration, supporting it to be efficient for
solving the large scale and sparse problems. Concretely, see Section
5 (Numerical exmaples). Moreover, it is different from
the following projection-type method for solving the
EHLCP($\mathrm{H},\mathrm{d}$) (\ref{3.9}):

\textbf{Method 3.3} \label{m33} (see Algorithm 4.1 \cite{Ahn}) For
$\mathrm{x}_{1}^{0}\in [0,\mathrm{b}]$, let
\begin{equation}\label{3.12}
\mathrm{x}_{1}^{k+1}=\eta
\mathbf{P}_{\mathrm{b}}[\mathrm{x}_{1}^{k}-\omega
\mathrm{E}(\mathrm{H}_{1}\mathrm{x}_{1}^{k}+\mathrm{q}+\mathrm{K}(\mathrm{x}_{1}^{k+1}-\mathrm{x}_{1}^{k}))]+(1-\eta)\mathrm{x}_{_{1}}^{k},
k=0, 1, \ldots,
\end{equation}
where $\mathbf{P}_{\mathrm{b}}$ is the projection operation on
$[0,\mathrm{b}]$, $0<\eta\leq1$ and $\omega>0$, $\mathrm{E}$ is any
positive diagonal matrix and $\mathrm{K}$ is either the strictly
lower or the strictly upper triangular part of $\mathrm{H}_{1}$.

Comparing Method 3.2 with Method 3.3, we know that
Method 3.2 can get the unique solution of the EHLCP($\mathrm{H},\mathrm{d}$) (\ref{3.9}) directly. 
In addition, Method 3.3 involves many parameters, it brings
inconvenience for computations.

\vspace{0.3cm}

The convergence theorem of Method 3.2 is given below. The proof is
similar to Theorem \ref{th33}, so we omit it.

\begin{theorem}\label{th37}
Let $\Omega=w\mathrm{I}$ with $w>0$. If
\begin{equation}\label{3.13}
\rho(|w^{-1}\mathrm{H}_{1}-\mathrm{I}|)<1,
\end{equation}
then the iterative sequence $\{\mathrm{y}^{k}\}^{+\infty}
_{k=0}\subset\mathbb{R}^{n}$ generated by Method 3.2 converges to
the unique solution $\mathrm{y}^{\ast}\in \mathbb{R}^{n}$ of the PLS
(\ref{3.10}) for an arbitrary initial vector
$\mathrm{y}^{0}\in\mathbb{R}^{n}$.
\end{theorem}
%

Clearly, we use the following simple condition
\begin{equation} \label{3.14}
\|w^{-1}\mathrm{H}_{1}-\mathrm{I}\|<1
\end{equation}
to replace  (\ref{3.13}) in order to guarantee convergence of Method
3.2.

Noting that the convergence conditions (\ref{3.13}) and (\ref{3.14})
do not demand that  $\mathrm{H}_{1}$ is a $\mathrm{Z}$-matrix or an
$\mathrm{H}_{+}$-matrix. In fact, taking
\[
\mathrm{H}_{1}=\left(\begin{array}{ccc}
1.5&1&1\\
1&1.5&1\\
1&1&1.5\\
\end{array}\right),
\]
it is easy to see that $\mathrm{H}_{1}$ is neither
$\mathrm{Z}$-matrix nor $\mathrm{H}_{+}$-matrix, but a $P$-matrix.
Set $w=5$, and thus we obtain
\[
\|w^{-1}\mathrm{H}_{1}-\mathrm{I}\|_{2}=0.90<1.
\]
In such case, we can not employ the modified Chandrasekaran's method in \cite{Pang79} and Method 3.3 
because the former  is a custom-tailored  direct method  for this
type of $\mathrm{Z}$-matrix and the latter one is a custom-tailored
iterative method for this type of $\mathrm{H}_{+}$-matrix. In
addition, this example also shows that the convergence conditions
(\ref{3.13}) and (\ref{3.14}) are not included each other because
$\rho(|w^{-1}\mathrm{H}_{1}-\mathrm{I}|)=1.1>1$.

A natural question is  how to choose the parameter $w$. In order to
answer this question, in the following we discuss it with two cases:
$\mathrm{H}_{1}$ is symmetric positive definite, and
$\mathrm{H}_{1}$ is a column sdd matrix with
$\mathrm{D}(\mathrm{H}_{1})=\tau \mathrm{I}$ and $\tau>0$, where
$\mathrm{D}(\mathrm{H}_{1})$ denotes the diagonal part of
$\mathrm{H}_{1}$, which frequently appears in many application
problems.

\begin{itemize}
\item $\mathrm{H}_{1}$ is symmetric positive definite. From   (\ref{3.14}), we have
\[
\|w^{-1}\mathrm{H}_{1}-\mathrm{I}\|_{2}=\max_{\delta_{i}\in
\delta(\mathrm{H}_{1})}\bigg|\frac{\delta_{i}-w}{w}\bigg|,
\]
where $\delta(\mathrm{H}_{1})$ is the spectrum of the matrix
$\mathrm{H}_{1}$. Clearly, we choose
$w>\frac{1}{2}\rho(\mathrm{H}_{1})$ to ensure the convergence of
Method 3.2. 
In practice, for saving the CPU time, we may choose
$w>\frac{1}{2}\|\mathrm{H}_{1}\|_{\infty}$ because of
$\|\mathrm{H}_{1}\|_{\infty}\geq\rho(\mathrm{H}_{1})$.

\item $\mathrm{H}_{1}$ is a column sdd matrix with
$\mathrm{D}(\mathrm{H}_{1})=\tau \mathrm{I}$ and $\tau>0$. For this
case, we choose $w=\tau$. In fact, it is easy to check that
\[
\|w^{-1}\mathrm{H}_{1}-\mathrm{I}\|_{1}<1.
\]
In addition, when the elements of $\mathrm{D}(\mathrm{H}_{1})$ are
positive, we take $\Omega=\mathrm{D}(\mathrm{H}_{1})$.
\end{itemize}

\section{The global error bound}
In this section, we focus on discussing the global error bounds of
the EHLCP($\mathrm{H},\mathrm{d}$). Assuming that
$\mathrm{w}^{\ast},\mathrm{x}^{\ast}_{1}$, \ldots,
$\mathrm{x}^{\ast}_{m}\in \mathbb{R}^{n}$ is a unique solution of
EHLCP($\mathrm{H},\mathrm{d}$). Let
\begin{equation} \label{4.1a} \tag{4.1a}
\begin{split}
\mathrm{r}(\mathrm{w},\mathrm{x}_{1}, \ldots,
\mathrm{x}_{m}):=\mathrm{q}+\sum_{i=1}^{m}
\mathrm{H}_{i}\mathrm{x}_{i}-\mathrm{Mw}
 \end{split}
 \end{equation}
 subject to
\addtocounter{equation}{1}
\begin{align}
& \mathrm{w},\mathrm{x}_{i}\geq0,i=1,\ldots,m, \label{4.1b} \tag{4.1b}\\
&\mathrm{w}^{T}\mathrm{x}_{1}=0,\label{4.1c} \tag{4.1c}\\
& \mathrm{x}_{i}\leq \mathrm{d}_{i}\ (\mathrm{d}_{i}>0),
(\mathrm{d}_{i}-\mathrm{x}_{i})^{T}\mathrm{x}_{i+1}=0,
i=1,2,\ldots,m-1. \label{4.1d} \tag{4.1d}
\end{align}
Clearly, $\mathrm{r}(\mathrm{w}^{\ast},\mathrm{x}^{\ast}_{1},
\ldots, \mathrm{x}^{\ast}_{m}) = 0$. If
$\mathrm{r}(\mathrm{w},\mathrm{x}_{1}, \ldots, \mathrm{x}_{m}) \neq
0$, then our purpose is to find the lower and upper bounds for
$\|\mathrm{r}(\mathrm{w},\mathrm{x}_{1}, \ldots, \mathrm{x}_{m})\|$.

Beyond their theoretical value for sensitivity analysis, the global
error bounds derived in this section have direct practical utility
for the iterative methods developed in Section 3. Specifically, the
residual $\mathrm{r}(\mathrm{w}(\mathrm{y}),
\mathrm{x}_1(\mathrm{y}), \ldots, \mathrm{x}_m(\mathrm{y}))$ defined
in (\ref{4.1a}) can be computed during the iteration of Methods 3.1
or 3.2. The bounds in Theorem 4.1, and their computable versions in
Theorems 4.2, 4.3 and 4.4, can then provide real-time estimates of
the distance to the true solution $\|\mathrm{y} - \mathrm{y}^*\|$,
serving as informative stopping criteria or solution quality
certificates. This establishes a practical connection between
algorithm design and error analysis, closing the theoretical loop
between Sections 3 and 4.

From Proposition \ref{pro21}, $\mathrm{r}(\mathrm{w},\mathrm{x}_{1},
\ldots, \mathrm{x}_{m})$ can be written as
\begin{align}
\mathrm{r}(\mathrm{w}(\mathrm{y}),\mathrm{x}_{1}(\mathrm{y}),
\ldots, \mathrm{x}_{m}(\mathrm{y}))=
&\sum_{i=1}^{m-1} \mathrm{H}_{i}\max\left \{ 0,\min\left \{ \mathrm{y}-\sum_{j=1}^{i-1} \mathrm{d}_{j},\mathrm{d}_{i} \right \} \right \} \label{4.2} \tag{4.2}\\
&+\mathrm{H}_{m}\max\left \{ 0,\mathrm{y}-\sum_{i=1}^{m-1}
\mathrm{d}_{i} \right \}-\mathrm{M}\max\{0,-\mathrm{y}\}+\mathrm{q},
\nonumber
\end{align}
where  $\mathrm{w}(\mathrm{y}),  \mathrm{x}_{1}(\mathrm{y}),\ldots,
\mathrm{x}_{m}(\mathrm{y})$ are defined as in (\ref{2.3}).

By Proposition \ref{pro21}, if
$\mathrm{w}^{\ast},\mathrm{x}^{\ast}_{1}, \ldots, x^{\ast}_{m}$ is a
unique solution of the EHLCP($\mathrm{H},\mathrm{d}$), then there
exists a corresponding vector $\mathrm{y}^{\ast}$ such that
(\ref{2.4}) holds, i.e., $\mathrm{r}(\mathrm{y}^{\ast})=0$, and
from Lemma \ref{y32}, 
there exists
$(\mathrm{D}_{0},\mathrm{D}_{1},\ldots,\mathrm{D}_{m})\in
\mathcal{D}$ such that 
\[
\mathrm{r}(\mathrm{w}(\mathrm{y}),\mathrm{x}_{1}(\mathrm{y}),
\ldots,
\mathrm{x}_{m}(\mathrm{y}))=(\mathrm{MD}_{0}+\mathrm{H}_{1}\mathrm{D}_{1}+\ldots+\mathrm{H}_{m}\mathrm{D}_{m})(\mathrm{y}-\mathrm{y}^{\ast}),
\]
which together with Lemma 1.1($ii$) gives the following theorem:

\begin{theorem}\label{th41} Let $\mathcal{H}$ have the column $\mathrm{W}$-property. Then
\begin{description}
  \item $(1)$ for any $(\mathrm{D}_{0},\mathrm{D}_{1},\ldots,\mathrm{D}_{m})\in \mathcal{D}$, the matrix $\mathrm{MD}_{0}+\mathrm{H}_{1}\mathrm{D}_{1}+\ldots+\mathrm{H}_{m}\mathrm{D}_{m}$ is nonsingular;
  \item $(2)$ for any $\mathrm{y}\in \mathbb{R}^{n}$ we have
\begin{equation}\label{4.3} \tag{4.3}
\underline{\alpha}^{-1}(\mathrm{H})\|\mathrm{r}(\mathrm{w}(\mathrm{y}),\mathrm{x}_{1}(\mathrm{y}),
\ldots, \mathrm{x}_{m}(\mathrm{y}))\|
\leq\|\mathrm{y}-\mathrm{y}^{\ast}\|\leq\overline{\alpha}(\mathrm{H})\|\mathrm{r}(\mathrm{w}(\mathrm{y}),\mathrm{x}_{1}(\mathrm{y}),
\ldots, \mathrm{x}_{m}(\mathrm{y}))\|,
\end{equation}
where
\begin{equation*}
\overline{\alpha}(\mathrm{H})=\max_{(\mathrm{D}_{0},\ldots,
\mathrm{D}_{m})\in
\mathcal{D}}\|(\mathrm{MD}_{0}+\mathrm{H}_{1}\mathrm{D}_{1}+\ldots+\mathrm{H}_{m}\mathrm{D}_{m})^{-1}\|
\end{equation*}
\mbox{and}
\begin{equation*}
\underline{\alpha}(\mathrm{H})=\max_{(\mathrm{D}_{0},\ldots,
\mathrm{D}_{m})\in
\mathcal{D}}\|\mathrm{MD}_{0}+\mathrm{H}_{1}\mathrm{D}_{1}+\ldots+\mathrm{H}_{m}\mathrm{D}_{m}\|.
\end{equation*}
\end{description}
\end{theorem}

The significance of Theorem \ref{th41} lies in its provision of a
unified framework for error analysis of the EHLCP. The bounds in
(\ref{4.3}) are expressed in terms of the transformed variable
$\mathrm{y}$, which is central to our iterative methods. The core
challenge reduces to estimating the norms of the parameter-dependent
matrices
$\mathrm{MD}_{0}+\mathrm{H}_{1}\mathrm{D}_{1}+\ldots+\mathrm{H}_{m}\mathrm{D}_{m}$
and their inverses. This framework  paves the way for deriving the
computable bounds presented in Theorems 4.2, 4.3 and 4.4 under
specific matrix structures.

\vspace{0.3cm}

For special cases we have

\begin{itemize}
\item For the HLCP($\mathrm{M}, \mathrm{H}_{1}, \mathrm{q}$) case, the bound (\ref{4.3}) reduces to the following one:

\begin{equation*}
\frac{1}{\underline{\beta}(\mathrm{H})}\|\mathrm{r}(\mathrm{w}(\mathrm{y}),\mathrm{x}_{1}(\mathrm{y}))\|
\leq\|\mathrm{y}-\mathrm{y}^{\ast}\|\leq\overline{\beta}(\mathrm{H})\|\mathrm{r}(\mathrm{w}(\mathrm{y}),\mathrm{x}_{1}(\mathrm{y}))\|,
\end{equation*}
where
\begin{equation*}
\overline{\beta}(\mathrm{H})=\max_{\mathrm{D}_{0}+\mathrm{D}_{1}=\mathrm{I}}\|(\mathrm{MD}_{0}+\mathrm{H}_{1}\mathrm{D}_{1})^{-1}\|\
\mbox{and}\
\underline{\beta}(\mathrm{H})=\max_{\mathrm{D}_{0}+\mathrm{D}_{1}=\mathrm{I}}\|\mathrm{MD}_{0}+\mathrm{H}_{1}\mathrm{D}_{1}\|,
\end{equation*}
and $\mathrm{D}_{0}, \mathrm{D}_{1}$ are nonnegative diagonal
matrices.

\item For the LCP($\mathrm{M},\mathrm{q}$) case, with the assumption of $\mathrm{M}$ being a $\mathrm{P}$-matrix,
the bound (\ref{4.3}) reduces to
\begin{equation*}
\frac{1}{\underline{\gamma}(\mathrm{H})}\|\mathrm{r}(\mathrm{w}(\mathrm{y}),
\mathrm{x}_{1}(\mathrm{y}))\|\leq\|\mathrm{y}-\mathrm{y}^{\ast}\|\leq\overline{\gamma}(\mathrm{H})\|\mathrm{r}(\mathrm{w}(\mathrm{y}),\mathrm{x}_{1}(\mathrm{y}))\|,
\end{equation*}
where $
\overline{\gamma}(\mathrm{H})=\max\|(\mathrm{MD}+\mathrm{I}-\mathrm{D})^{-1}\|\
\mbox{and}\
\underline{\gamma}(\mathrm{H})=\max\|\mathrm{MD}+\mathrm{I}-\mathrm{D}\|,
$
and $\mathrm{D}=\mbox{diag}(d)$ with $d\in [0,1]^{n}$ is nonnegative
diagonal matrix.

\end{itemize}


It is difficult to directly compute \(\overline{\alpha}(\mathrm{H})\) in general. Similar to Theorems 2 and 3 in [20], we can present some computable error bounds for the special block matrix \(\mathrm{H}\): (I) \(\mathrm{H} = (\mathrm{M}, \mathrm{H}_1, \dots, \mathrm{H}_m)\) with each matrix having positive diagonal part; (II) \(\mathrm{H} = (\mathrm{M}, \mathrm{H}_1, \dots, \mathrm{H}_m)\) with each matrix being a column sdd matrix. That is, we present Theorems 4.2 and 4.3 for these two cases, respectively. Their proofs follow directly from the analogous results in [20] and are thus omitted.

\begin{theorem}\label{th42}
Let  \(\mathrm{M} = \Lambda_0 - \mathrm{C}_0\), \(\mathrm{H}_i = \Lambda_i - \mathrm{C}_i\) be the splitting of \(\mathrm{M}\) and \(\mathrm{H}_i\), respectively, where \(\Lambda_i\) is the diagonal part of \(\mathrm{M}\) and \(\mathrm{H}_i\), \(i = 0, 1, \dots, m\). If \(\Lambda_i > 0\) and
\begin{equation} \label{4.4} \tag{4.4}
\rho \left( \max_{0 \le i \le m} \{ |\mathrm{C}_i| \Lambda_i^{-1} \} \right) < 1,
\end{equation}
then  \(\mathcal{H}\) has the column W-property, and
\begin{equation}\label{4.5} \tag{4.5}
\overline{\alpha}(\mathrm{H}) \le \left\|\max_{0 \le i \le m} \{ \Lambda_i^{-1} \} \left( \mathrm{I} - \max_{0 \le i \le m} \{ |\mathrm{C}_i|\Lambda_i^{-1} \} \right)^{-1} \right\|.
\end{equation}
\end{theorem}

\begin{theorem} \label{th43}
Let the block matrix $\mathrm{H}=(\mathrm{M},\mathrm{H}_{1},\ldots,
\mathrm{H}_{m})$, where $\mathrm{M}^{T}$ and $\mathrm{H}^{T}_{i}$
are the row sdd matrices with the $\ell$-th diagonal entry having
the same sign, $\ell = 1, \ldots, n$. Then $\mathcal{H}$ has the
column $\mathrm{W}$-property and
\begin{equation}\label{4.6} \tag{4.6}
\overline{\alpha}_{1}(\mathrm{H})\leq\frac{1}{\min_{i\in
\mathcal{N}}\min((\langle \mathrm{M}^{T}\rangle
\mathrm{e})_{i},(\langle \mathrm{H}^{T}_{1}\rangle
\mathrm{e})_{i},\ldots, (\langle \mathrm{H}^{T}_{m}\rangle
\mathrm{e})_{i})},
\end{equation}
where $\overline{\alpha}_{1}(\mathrm{H})$ denotes the 1-norm form of
$\overline{\alpha}(\mathrm{H})$,  i.e.,
\[
\bar{\alpha}_{1}(\mathrm{H}) :=
\max_{(\mathrm{D}_{0},\dots,\mathrm{D}_{m})\in\mathcal{D}} \|
(\mathrm{M}\mathrm{D}_{0}+\mathrm{H}_{1}\mathrm{D}_{1}+\dots+\mathrm{H}_{m}\mathrm{D}_{m})^{-1}
\|_{1},
\]
and $\| \cdot \|_{1}$ for a matrix is the maximum absolute column
sum.
\end{theorem}

\begin{remark}\label{rem:4.1}
It should be pointed out that the conditions in Theorems 4.2 and 4.3 are not mutually inclusive. For instance, taking the block matrix \(\mathrm{H} = (\mathrm{M}, \mathrm{H}_1)\), where
\[
\mathrm{M} = \begin{pmatrix} 1 & 0 \\ -1 & 1 \end{pmatrix}, \quad \mathrm{H}_1 = \begin{pmatrix} 1 & 0 \\ 2 & 1 \end{pmatrix}.
\]
A simple computing gives
\[
\rho(\max\{|C_0|\Lambda_0^{-1},|C_1|\Lambda_1^{-1}\}) = 0 < 1.
\]
Then \(\mathrm{H}\) satisfies the condition (4.4) in Theorem 4.2. It is easy to see that \(\mathrm{M}\) and \(\mathrm{H}_1\) are not column sdd matrices, from which \(\mathrm{H}\) does not satisfy the condition in Theorem 4.3. Now, we take \(\mathrm{H} = (\mathrm{M}, \mathrm{H}_1)\), where
\[
\mathrm{M} = \begin{pmatrix} 2 & 0 & 0 \\ 1 & 2 & 1 \\ 0 & 1 & 2 \end{pmatrix}, \quad \mathrm{H}_1 = \begin{pmatrix} 2 & 1 & 1 \\ 0 & 2 & 0 \\ 1 & 0 & 2 \end{pmatrix}.
\]
Then both \(\mathrm{M}\) and \(\mathrm{H}_1\) are column sdd matrices. A simple computing gives
\[
\rho(\max\{|C_0|\Lambda_0^{-1},|C_1|\Lambda_1^{-1}\}) = 1.
\]
This shows that \(\mathrm{H}\) does not satisfy the condition (4.4) in Theorem 4.2.
\end{remark}

\section{Numerical examples}
In this section,  we will show the feasibility and efficiency of the
proposed error bounds and the convergence behavior of Method 3.2.
All the numerical experiments are computed by using Matlab R2024a on
the Lenvo PC (13th Gen Intel(R) Core(TM) i7-13700, 2.10GHz, 16.00 GB
of RAM). The test problems are chosen to be large-scale and sparse,
reflecting the intended application domain of the proposed methods.
The algorithms' performance in terms of iteration count and CPU time
will demonstrate their efficiency and suitability for such problems.

\subsection{The numerical test for error upper bounds}
In this subsection, our task is to  confirm the feasibility and
efficiency of the upper bounds of Theorems \ref{th42} and
\ref{th43}.

\textbf{Example 5.1} (\cite{Mezzadri2}) Consider the
HLCP($\mathrm{M}, \mathrm{H}_{1}, \mathrm{q}$), in which
$\mathrm{M}=\bar{\mathrm{A}}+\mu \mathrm{I} \in \mathbb{R}^{n\times
n}$, $\mathrm{H}_{1}=\bar{\mathrm{B}}+\nu \mathrm{I} \in
\mathbb{R}^{n\times n}$, where
$\bar{\mathrm{A}}=\mbox{blktridiag}(-\mathrm{I},\mathrm{T},-\mathrm{I})$,
$\bar{\mathrm{B}}=\mathrm{I}\otimes \mathrm{T}$,
$\mathrm{T}=\mbox{tridiag}(-1,4,-1)\in \mathbb{R}^{m\times m}$ (a
sparse tridiagonal matrix), $n=m^{2}$ and $\mu\geq0 $, $\nu\geq0$.

This HLCP($\mathrm{M}, \mathrm{H}_{1}, \mathrm{q}$) can be derived
from the discretization of certain differential equation under the
complementarity conditions. For instance, under some 2D domain of
axes $(x, y)$, one considers
\begin{equation*}
\Delta u(x,y)+\frac{\partial^{2}v(x,y)}{\partial x^{2}}+\mu
u(x,y)+\nu v(x,y)=q(x,y), u\geq0, v\geq0, u^{T}v=0
\end{equation*}
with suitable boundary conditions, where $u(x, y), v(x, y)$ and
$q(x, y)$ are three 2-D mappings. By making use of the five-point
difference scheme for the approximation of the Laplacian and central
finite difference scheme for the approximation of the second-order
partial derivative, we can get the desired HLCP($\mathrm{M},
\mathrm{H}_{1}, \mathrm{q}$).

In our computations, we set the vector $\mathrm{q}$ such that
$\mathrm{q}=\mathrm{Mw}^{\ast}-\mathrm{Hx}_{1}^{\ast}$ with
\[
\mathrm{w}^{\ast}=(0.1,0,0.1,0 \ldots,0,0.1,\ldots)^{T}\in
\mathbb{R}^{n},\mathrm{x}_{1}^{\ast}=(0, 0.1, 0,0.1
\ldots,0,0.1,\ldots)^{T}\in \mathbb{R}^{n}.
\]
where $\mathrm{w}^{\ast}, \mathrm{x}_{1}^{\ast}$ is the solution of
the HLCP($\mathrm{M}, \mathrm{H}_{1}, \mathrm{q}$). By Theorem
\ref{th32},
\[
\mathrm{y}^{\ast}=\mathrm{x}_{1}^{\ast}-\mathrm{w}^{\ast}=(-0.1,0.1,\ldots,-0.1,0.1,\ldots)^{T}\in
\mathbb{R}^{n}
\]
is the unique solution of the corresponding PLS (\ref{2.4}).

For convenience, we take $\mu=\nu>0$. In this setting, $\mathrm{M}$
and $\mathrm{H}_{1}$ are two symmetric column sdd matrices. A simple
computation reveals that
\begin{align*}
\bar{\eta}=\|\max\{\mathrm{\Lambda}^{-1}_{0},\mathrm{\Lambda}^{-1}_{1}\}(\mathrm{I}-
\max\{\mathrm{|\mathrm{C}_{0}|\Lambda}^{-1}_{0},|\mathrm{C}_{1}|\mathrm{\Lambda}^{-1}_{1}\}
)^{-1}\|_{\infty}
=\|(\mathrm{\Lambda}_{0}-|\mathrm{C}_{0}|)^{-1}\|_{\infty}=\|\langle
\mathrm{M}\rangle^{-1}\|_{\infty}
\end{align*}
and
\begin{align*}
\bar{\tau}=\frac{1}{\min_{i\in \mathcal{N}}\min((\langle
\mathrm{M}^{T}\rangle e)_{i},(\langle \mathrm{H}^{T}_{1}\rangle
e)_{i})} =\frac{1}{\min_{i\in \mathcal{N}}\min((\langle
\mathrm{M}\rangle e)_{i},(\langle \mathrm{H}_{1}\rangle e)_{i})}.
\end{align*}
This implies that $\bar{\eta}\leq\bar{\tau}$. Let
\[
\mathrm{r}_{\infty}=\|\mathrm{y}-\mathrm{y}^{\ast}\|_{\infty},\eta_{\infty}=\bar{\eta}\|\mathrm{r}(\mathrm{w}(\mathrm{y}),\mathrm{x}_{1}(\mathrm{y}))\|_{\infty},
\tau_{\infty}=\bar{\tau}\|\mathrm{r}(\mathrm{w}(\mathrm{y}),\mathrm{x}_{1}(\mathrm{y}))\|_{\infty}.
\]
Clearly, $\eta_{\infty}\leq \tau_{\infty}$. Moreover, it is easy to find that  the larger $\mu$ is, the smaller $\eta_{\infty}$ and $\tau_{\infty}$ are. 
To further verify it, for the vector
\[
\mathrm{y}=(-0.15,0.056,\ldots,-0.15,0.056,\ldots)^{T}\in
\mathbb{R}^{n},
\]
we report the numerical results in the following table.

\begin{table}[!htb] \centering
\begin{tabular}
{p{10pt}p{35pt}p{35pt}p{35pt}p{35pt}p{35pt}p{35pt}} \hline
$\mu$&4&6&8&10&12&14
\\\hline
$r_{\infty}$    &  0.05000  &  0.05000 &  0.05000  &  0.05000 &  0.05000 &  0.05000\\
$\eta_{\infty}$ &  0.07650 &   0.06767 &  0.06325 &0.06060    & 0.05883 &0.05757\\
$\tau_{\infty}$ &  0.07650 &    0.06767 & 0.06325 &0.06060  &0.05883 &0.05757\\
\hline
\end{tabular}
\\ \caption{Numerical results of Example 5.1 with $n=10000$.}
\end{table}

In addition, from the numerical results in Tab. 5.1, as $\mu$
increases, $\eta_{\infty}$ and $\tau_{\infty}$ decrease and
gradually approach to $r_{\infty}$.

\begin{figure}
\subfigure[$\mu=5$]{
\begin{minipage}[b]{0.23\textwidth}
\centering
\includegraphics[width=1.5in]{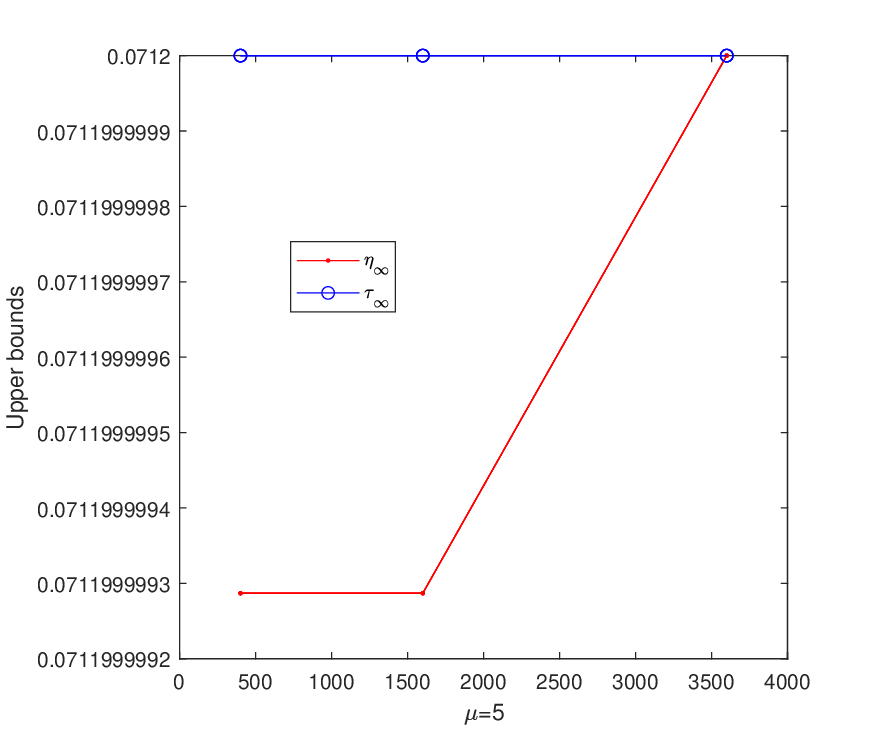}
\end{minipage}}%
 \subfigure[$\mu=7$]{
\begin{minipage}[b]{0.23\textwidth}
\centering
\includegraphics[width=1.5in]{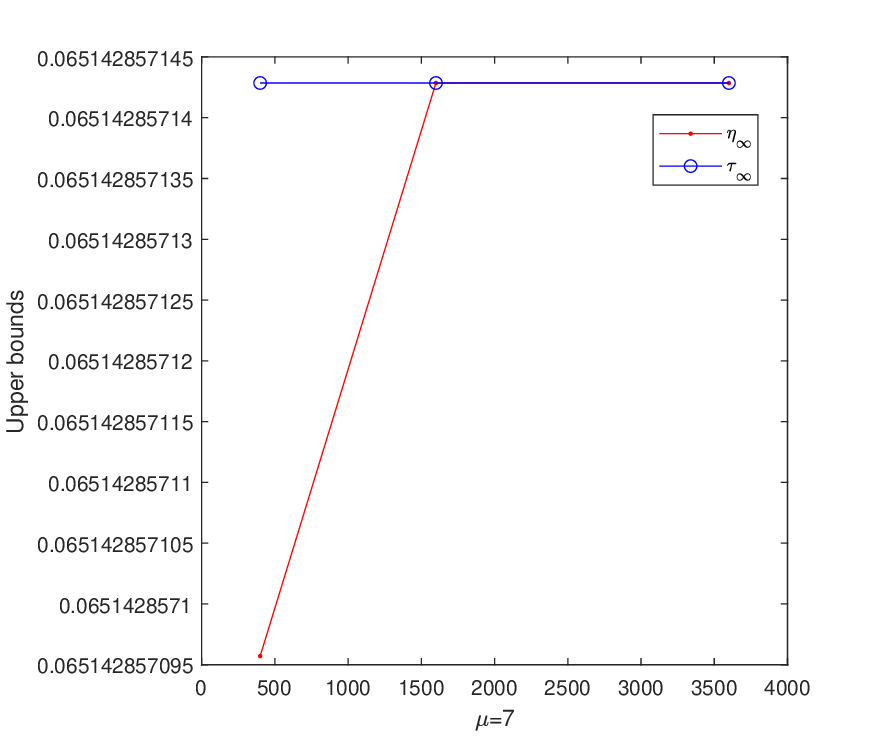}
\end{minipage}}
\subfigure[$\mu=9$]{
\begin{minipage}[b]{0.23\textwidth}
\centering
\includegraphics[width=1.5in]{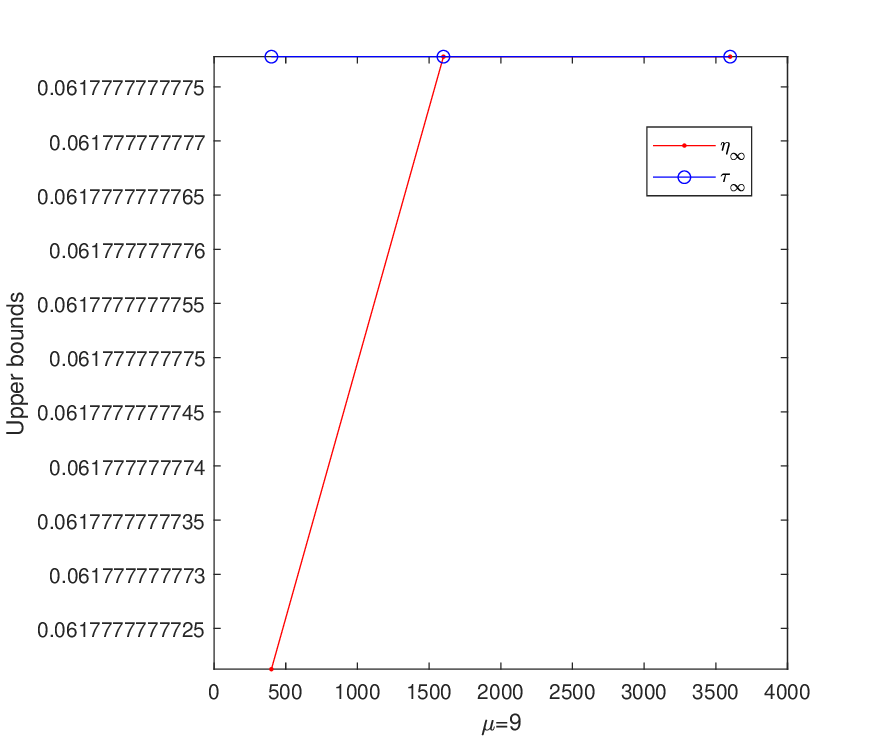}
\end{minipage}}
\subfigure[$\mu$]{
\begin{minipage}[b]{0.23\textwidth}
\centering
\includegraphics[width=1.45in]{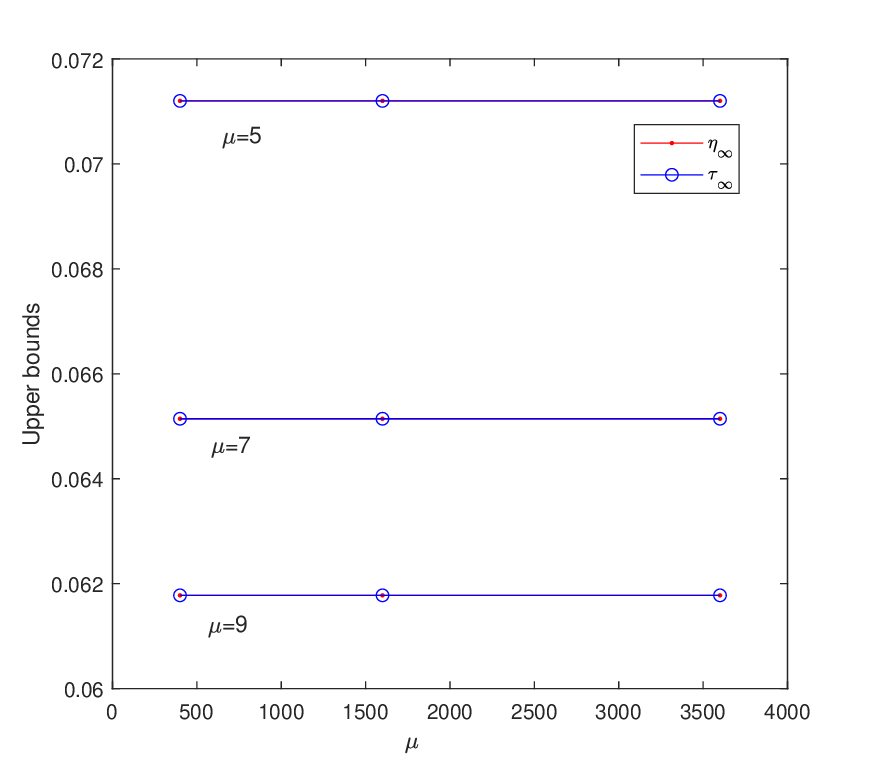}
\end{minipage}}
\caption{The upper bounds of $\mu$ in Table 2.} \label{f2}
\end{figure}

Next, we report the numerical results of the error bounds for
testing the effect on mesh size $n$ in Tab. 5.2, in which we can see
that for a given $\mu$, with the mesh size $n$ increasing,
$\eta_{\infty}$ increases, but quite slowly, and $\tau_{\infty}$ has
no change, also see (a), (b) and (c) in Fig. \ref{f2}. In Fig.
\ref{f2} (d) we plot the effect of coupling (a), (b) and (c). Based
on these numerical results, both $\tau$ and $\mu$  are insensitive
to change in the mesh sizes for a given $\mathrm{y}$.

\begin{table}[htbp]
\centering
\begin{tabular}
{|p{10pt}|p{10pt}|p{100pt}|p{100pt}|p{100pt}|} \hline $\mu$&
$n$&400&1600&3600\\ \hline \raisebox{-1.50ex}[0cm][0cm]{5}&
$\eta_{\infty}$
&0.071199999286907&0.071199999286907&0.071200000000000
\\
\cline{2-5}
&$\tau_{\infty}$&0.071200000000000&0.071200000000000&0.071200000000000
\\
\cline{1-5}  \raisebox{-1.50ex}[0cm][0cm]{7}&$\eta_{\infty}$&  0.065142857095714&0.065142857142857&0.065142857142857\\
\cline{2-5} &$\tau_{\infty}$& 0.065142857142857&0.065142857142857&0.065142857142857\\
\cline{1-5}  \raisebox{-1.50ex}[0cm][0cm]{9}&$\eta_{\infty}$& 0.061777777772124& 0.061777777777778&0.061777777777778\\
\cline{2-5} &$\tau_{\infty}$&0.061777777777778& 0.061777777777778&0.061777777777778\\
\hline
\end{tabular}
\caption{Numerical results of Example 5.1 for the different sizes
$n$ and values of $\mu$.}
\end{table}

Tabs. 5.1 and 5.2 imply that the upper bounds given by Theorems \ref{th42} and \ref{th43} are feasible and effective. 

\bigskip

\textbf{Example 5.2}\label{ex52} (\cite{Pang79,Ahn}) Consider the
EHLCP($\mathrm{H},\mathrm{d}$) (\ref{3.9}),
which can arise from multicommodity market equilibrium problem with
institutional price controls imposed.

Here, according to the approach in \cite{Ahn}, for Example 5.2,  we
examine the upper bounds of Theorems \ref{th42} and \ref{th43} for
$\mathrm{H}_{1}=\mbox{tridiag}(1,4,-2)\in \mathbb{R}^{n\times n}$ (a
sparse tridiagonal matrix). For convenience, we choose
$\mathrm{b}=0.1\mathrm{e}$, and set the vector $\mathrm{q}$ such
that
$\mathrm{q}=\mathrm{w}^{\ast}-\mathrm{H}\mathrm{x}_{1}^{\ast}-\mathrm{x}_{2}$
with
\[
\mathrm{w}^{\ast}=(0.2,0,0.2 \ldots,0.2,0,\ldots)^{T}\in
\mathbb{R}^{n},\mathrm{x}_{1}^{\ast}=\mathrm{x}_{2}^{\ast}=(0, 0.1,0
\ldots,0,0.1,\ldots)^{T}\in \mathbb{R}^{n},
\]
where $\mathrm{w}^{\ast}, \mathrm{x}_{1}^{\ast},
\mathrm{x}_{2}^{\ast}$ is the solution of (\ref{3.9}). Then by
Theorem \ref{th32},
\[
\mathrm{y}^{\ast}=\mathrm{x}_{1}^{\ast}+\mathrm{x}_{2}^{\ast}-\mathrm{w}^{\ast}=(-0.2,0.2,\ldots,-0.2,0.2,\ldots)^{T}\in
\mathbb{R}^{n}
\]
is the unique solution of the corresponding PLS (\ref{2.4}).

Similar to Example 5.1, by a simple computation, we get
\[
\hat{\eta}=\|(\mathrm{I}-\mathrm{|\mathrm{C}_{1}|\Lambda}^{-1}_{1})^{-1}\|_{\infty},
\hat{\tau} =\frac{1}{\min_{i\in \mathcal{N}}\min(1,(\langle
\mathrm{H}^{T}_{1}\rangle e)_{i})}.
\]
Since $\mathrm{H}_{1}$ is unsymmetrical, we consider $r_{\infty}$
and $r_{1}$ cases for Theorems \ref{th42} and \ref{th43},
respectively. Let  $r_{\infty}=\|\mathrm{y}-\mathrm{y}^{\ast}\|_{\infty}$, $r_{1}=\|\mathrm{y}-\mathrm{y}^{\ast}\|_{1},$
and
\[
\eta_{\infty}=\hat{\eta}\|\mathrm{r}(\mathrm{w}(\mathrm{y}),\mathrm{x}_{1}(\mathrm{y}),\mathrm{x}_{2}(\mathrm{y}))\|_{\infty},
\tau_{1}=\hat{\tau}\|\mathrm{r}(\mathrm{w}(\mathrm{y}),\mathrm{x}_{1}(\mathrm{y}),\mathrm{x}_{2}(\mathrm{y}))\|_{1}.
\]
Then for the vector
\[
\mathrm{y}=(-0.1,0.1,\ldots,-0.1,0.1,\ldots)^{T}\in \mathbb{R}^{n},
\]
we list some numerical results in Tabs. 5.3 and 5.4.

\begin{table}[!htb] \centering
\begin{tabular}
{p{10pt}p{35pt}p{35pt}p{35pt}p{35pt}p{35pt}p{35pt}} \hline
$n$&30&60&90&120
\\\hline
$r_{1}$    &  3  &  6  &  9  &  12 \\
$\tau_{1}$ & 3 &   6&  9 &12  \\
\hline
\end{tabular}
\\ \caption{Numerical results of Example 5.2 with $r_{1}$ and $\tau_{1}$.}
\end{table}

\begin{table}[!htb] \centering
\begin{tabular}
{p{10pt}p{35pt}p{35pt}p{35pt}p{35pt}p{35pt}p{35pt}} \hline
$n$&30&60&90&120
\\\hline
$r_{\infty}$    &  0.1  &  0.1  &  0.1  &  0.1 \\
$\eta_{\infty}$ &  0.4 &   0.4&  0.4 &0.4  \\
\hline
\end{tabular}
\\ \caption{Numerical results of Example 5.2 with $r_{\infty}$ and $\eta_{\infty}$.}
\end{table}

It can be seen from Tab. 5.3 that for
$\mathrm{y}=(-0.1,0.1,\ldots,-0.1,0.1,\ldots)^{T}\in
\mathbb{R}^{n}$, with the mesh size $n$ increasing, $r_{1}$ and
$\tau_{1}$ increase. A more important fact is that $r_{1}$ and
$\tau_{1}$ are the same. This implies that using Theorem \ref{th43}
to estimate the global error bounds is sharper.

From Tab. 5.4, we can see that for
$\mathrm{y}=(-0.1,0.1,\ldots,-0.1,0.1,\ldots)^{T}\in
\mathbb{R}^{n}$, with the mesh size $n$ increasing, $r_{\infty}$ and
$\eta_{\infty}$ have no change, which shows that they are
insensitive to the mesh sizes. Clearly, $r_{\infty}<\eta_{\infty}$.
This shows that the upper bound given by Theorem \ref{th42} is
feasible
and effective under some suitable conditions. 

\bigskip

A simple example is given below to illustrate the superiority of the
global error bounds in Theorem \ref{th42}.

\textbf{Example 5.3} For $\alpha\in [1,+\infty)$, we consider
$\mathrm{H}=(\mathrm{M},\mathrm{H}_{1})$, where
\[
\mathrm{M}=\left(\begin{array}{ccc}
1&0\\
\alpha&1\\
\end{array}\right),\ \mathrm{H}_{1}=\left(\begin{array}{ccc}
1&0\\
\alpha^{2}&1\\
\end{array}\right).
\]

A simple computation gives
\[
\rho(\max\{\mathrm{|\mathrm{C}_{0}|\Lambda}^{-1}_{0},|\mathrm{C}_{1}|\mathrm{\Lambda}^{-1}_{1}\})=0<1.
\]
This shows that  $\mathrm{H}$ satisfies the condition in Theorem
\ref{th42}, and has the column $\mathrm{W}$-property. Setting
\[
\mathrm{D}_{0}=\left(\begin{array}{ccc}
1-d^{1}_{0}&0\\
0&1-d^{2}_{0}\\
\end{array}\right).
\]
Then we have
\begin{align*}
\overline{\beta}_{\infty}(\mathrm{H})&=\max_{\mathrm{D}_{0}+\mathrm{D}_{1}=I}\|(\mathrm{MD}_{0}+\mathrm{H}_{1}\mathrm{D}_{1})^{-1}\|_{\infty}\\
&=\max_{d^{1}_{0}\in [0,1]}(1+\alpha
(1-d^{1}_{0})+\alpha^{2}d^{1}_{0})\\
&=\|\max_{0\leq i\leq 1}\{\mathrm{\Lambda}^{-1}_{i}\}(\mathrm{I}-\max_{0\leq i\leq 1}\{\mathrm{|\mathrm{C}_{i}|\Lambda}^{-1}_{i}\})^{-1}\|_{\infty}\\
&=1+\alpha^{2}.
\end{align*}

Let $\alpha=1$, and $\mathrm{q}=(1,0)^{T}$. Then
$\mathrm{w}^{\ast}=(1,0)^{T}, \mathrm{x}_{1}^{\ast}=(0,1)^{T}$ is
the unique solution of the corresponding HLCP($\mathrm{M},
\mathrm{H}_{1}, \mathrm{q}$), and $\mathrm{y}^{\ast}=(-1,1)^{T}$
from Theorem \ref{th32} is the unique solution of the corresponding
PLS (\ref{2.4}). For the vector $\mathrm{y}=(3,-7)^{T}$, we get
\[
\|\mathrm{y}-\mathrm{y}^{\ast}\|_{\infty}=8,
\overline{\beta}_{\infty}(\mathrm{H})=2,
\|\mathrm{r}(\mathrm{w}(\mathrm{y}),\mathrm{x}_{1}(\mathrm{y}))\|_{\infty}=4.
\]
Clearly,
\[
\|\mathrm{y}-\mathrm{y}^{\ast}\|_{\infty}=\overline{\beta}_{\infty}(\mathrm{H})\|\mathrm{r}(\mathrm{w}(\mathrm{y}),\mathrm{x}_{1}(\mathrm{y}))\|_{\infty}.
\]
This shows that the bound in (\ref{4.3}) can be achieved.

\subsection{Test the proposed method}
In this subsection, we give some numerical examples to show the convergence behavior of Method 3.2. 
The termination criterion is set by
$\|\mathrm{y}^{k+1}-\mathrm{y}^{k}\|_{\infty}<10^{-6} $ for Method
3.2, and
$\|\mathrm{x}_{1}^{k+1}-\mathrm{x}_{1}^{k}\|_{\infty}<10^{-6} $ for
Method 3.3. All initial vectors are taken to zero. Set $\eta=0.5$,
$\omega=0.25$, $\mathrm{E}=\mathrm{I}$ and
$\mathrm{K}=-\mathrm{L}(\mathrm{H}_{1})$ for Method 3.3,
$\mathrm{L}(\mathrm{H}_{1})$ is the strictly lower triangular part
of $\mathrm{H}_{1}$. In the following tables, `IT' and `CPU' in
order denote the iterative step and CPU time (second) on a selection
of test problems. By `M2' and `M3' we denote Method 3.2 and Method
3.3, respectively.

\textbf{Example 5.4} Example 5.2 is used as the first example to
test the proposed algorithm. Since $\mathrm{H}_{1}$ is a column and
row sdd matrix, we set $w=4$ for Method 3.2.
In Tab. 5.5 we report numerical results for Example 5.4.

\begin{table}[!htb] \centering
\begin{tabular}
{p{20pt}p{35pt}p{35pt}p{35pt}p{35pt}p{35pt}p{35pt}} \hline
&$n$&5000&10000&15000&20000
\\\hline
M2 &  IT &  3  &  3 &  3  &  3 \\
   & CPU &  0.0011  & 0.0014   &  0.0015  & 0.0024  \\
M3 & IT& 16 &   16&  16 &16 \\
   & CPU & 2.7589  &   10.1852  &  22.4088  &  39.0979 \\
\hline
\end{tabular}
\\ \caption{Numerical results of Example 5.4.}
\end{table}

\begin{figure}
\subfigure[$n=5000$]{
\begin{minipage}[b]{0.23\textwidth}
\centering
\includegraphics[width=1.5in]{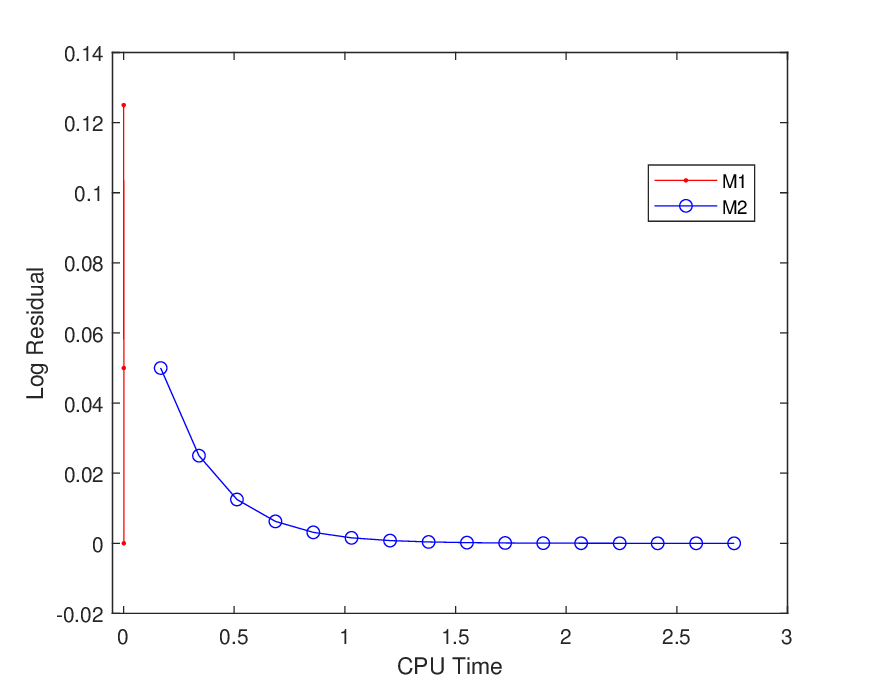}
\end{minipage}}%
 \subfigure[$n=1000$]{
\begin{minipage}[b]{0.23\textwidth}
\centering
\includegraphics[width=1.5in]{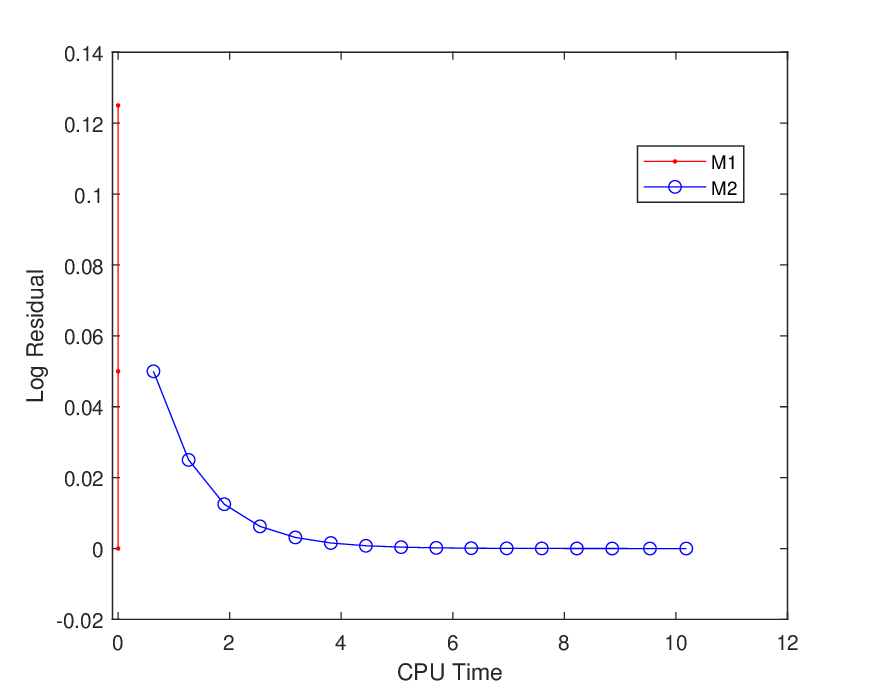}
\end{minipage}}
\subfigure[$n=15000$]{
\begin{minipage}[b]{0.23\textwidth}
\centering
\includegraphics[width=1.5in]{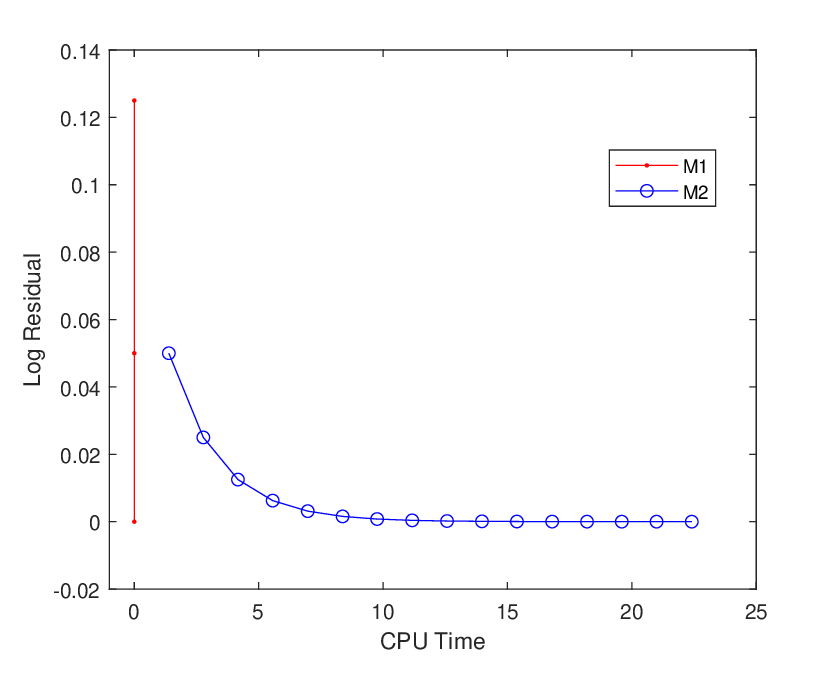}
\end{minipage}}
\subfigure[$n=20000$]{
\begin{minipage}[b]{0.23\textwidth}
\centering
\includegraphics[width=1.5in]{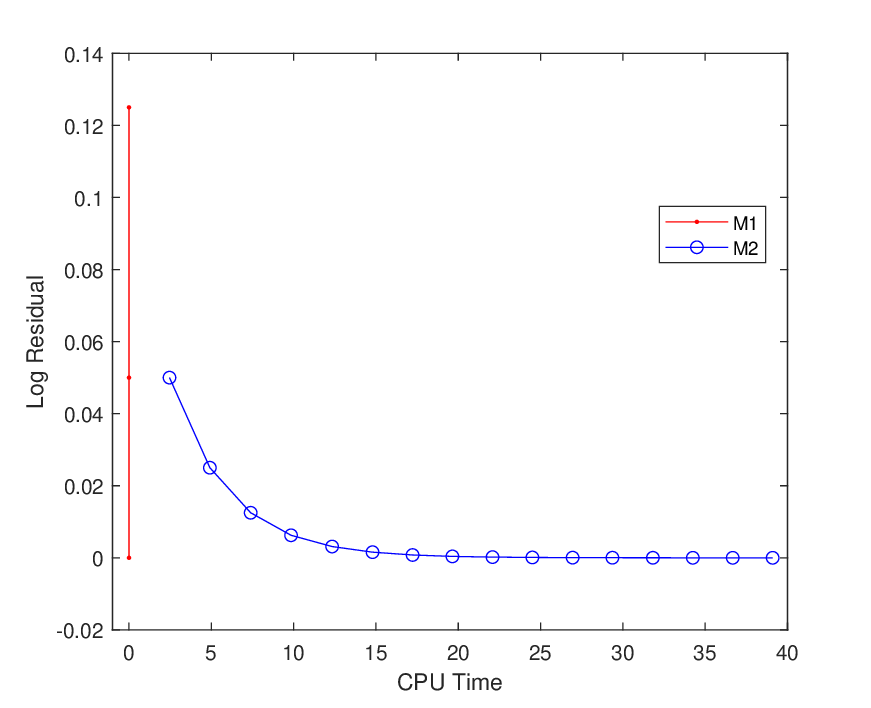}
\end{minipage}}
\caption{Errors and CPU times of Table 5.}\label{f3}
\end{figure}

\textbf{Example 5.5} \label{ex55} (\cite{Yuan12}) Consider the bilateral obstacle problem: 
\[
E(u)=\min_{v\in K}\frac{1}{2}\int_{\Omega} |\nabla
v|^{2}dxdy-\int_{\Omega}qvdxdy,
\]
where $\Omega\in \mathbb{R}^{2}$ is a bounded domain with Lipschitz
boundary  $\partial\Omega$, $K$ is a closed convex subset of
$H^{1}_{0}(\Omega)$  and
\[
K = \{v \in H^{1}_{0}(\Omega), \phi\leq v \leq \psi, \mbox{a.e.} \ \
\mbox{in} \ \Omega\},
\]
where the obstacle functions $\phi$ and $\psi$ satisfy $\phi\leq
\psi \in H_{1}(\Omega)\cap C(\Omega)$, and $\phi\leq0$ and
$\psi\geq0$ on $\partial\Omega$, $q\in L^{2}(\Omega)$, also see
\cite{Cottle78} for more details. Following the approach  in
\cite{Yuan12}, with only slight and technical modifications, we
obtain the EHLCP($\mathrm{H},\mathrm{d}$) (\ref{3.9}) with
$\mathrm{H}_{1}=\mbox{blktridiag}(-\mathrm{I},\mathrm{T},-\mathrm{I})\in
\mathbb{R}^{n\times n}$, $\mathrm{T}=\mbox{tridiag}(-1,4,-1)\in
\mathbb{R}^{m\times m}$, $n=m^{2}$ and $\mathrm{b}=0.1\mathrm{e}$.
We adjust the vector $\mathrm{q}$ such that the solution of Example
5.5 is the same as that of Example 5.2.

\begin{table}[!htb] \centering
\begin{tabular}
{p{20pt}p{35pt}p{35pt}p{35pt}p{35pt}p{35pt}p{35pt}} \hline
&$n$&6400&10000&16900&22500
\\\hline
M2 &  IT &  5  &  5 &  5  &  5 \\
   & CPU &  0.0026  & 0.0034   &  0.0042  & 0.0060  \\
M3 & IT& 16 &   16&  16 &16  \\
   & CPU &  5.0362  &   12.2314  &  34.2181  &  62.5915 \\
\hline
\end{tabular}
\\ \caption{Numerical results of Example 5.5.}
\end{table}

\begin{figure}
\subfigure[$n=6400$]{
\begin{minipage}[b]{0.23\textwidth}
\centering
\includegraphics[width=1.5in]{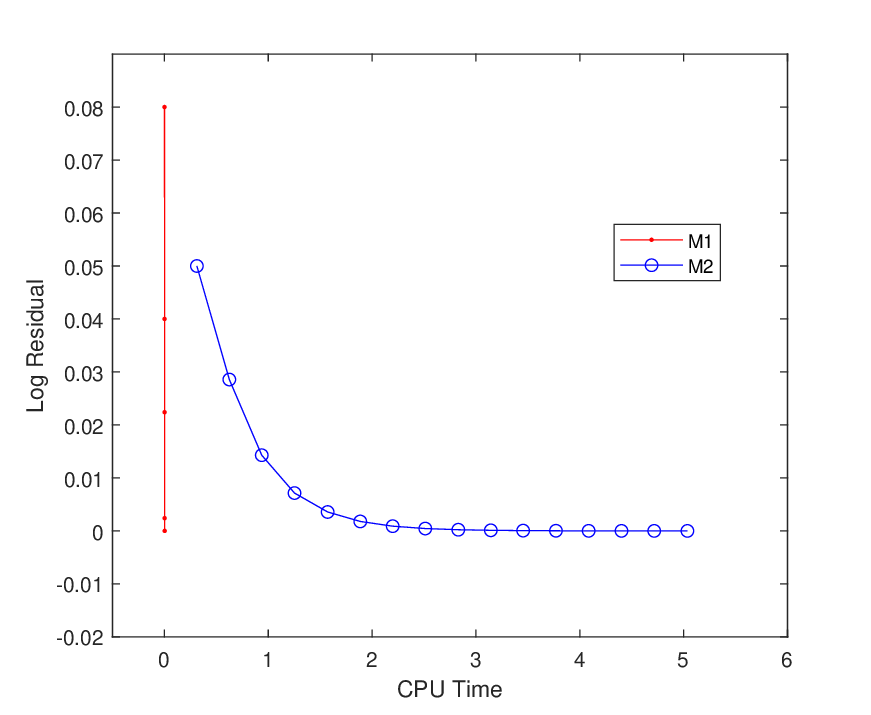}
\end{minipage}}%
 \subfigure[$n=10000$]{
\begin{minipage}[b]{0.23\textwidth}
\centering
\includegraphics[width=1.5in]{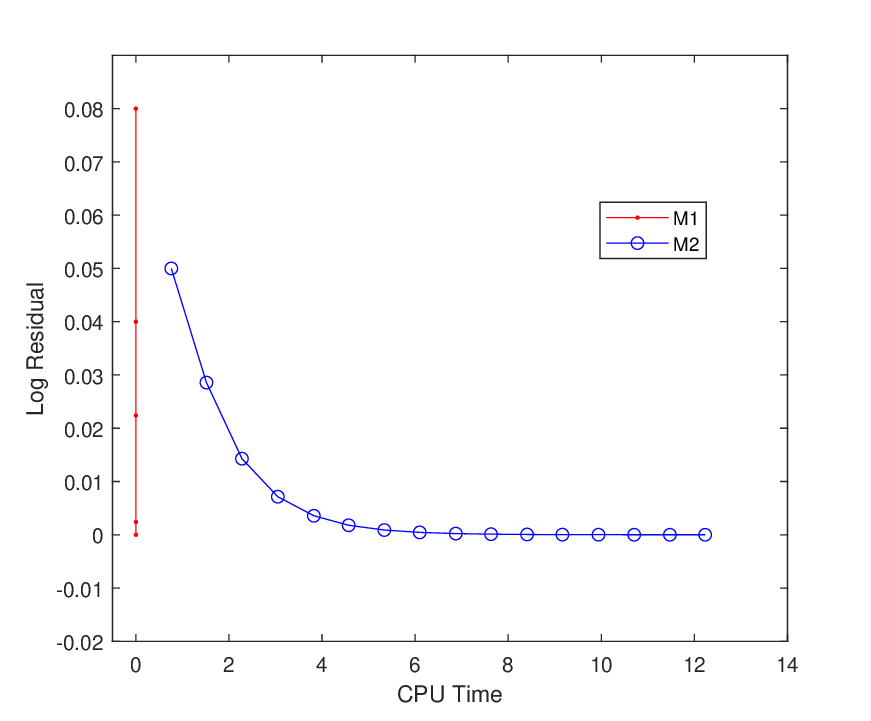}
\end{minipage}}
\subfigure[$n=16900$]{
\begin{minipage}[b]{0.23\textwidth}
\centering
\includegraphics[width=1.5in]{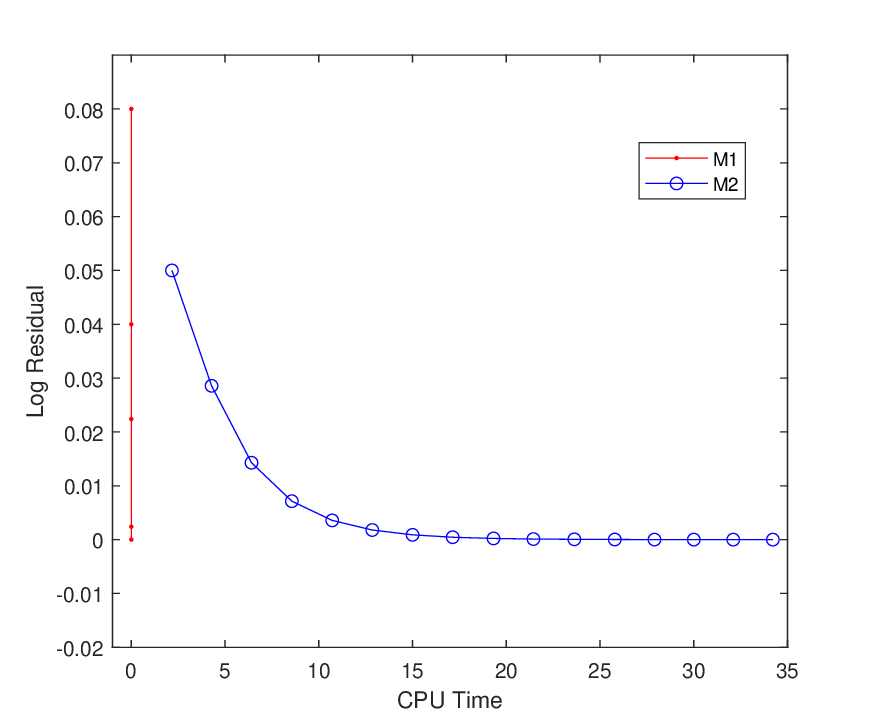}
\end{minipage}}
\subfigure[$n=22500$]{
\begin{minipage}[b]{0.23\textwidth}
\centering
\includegraphics[width=1.5in]{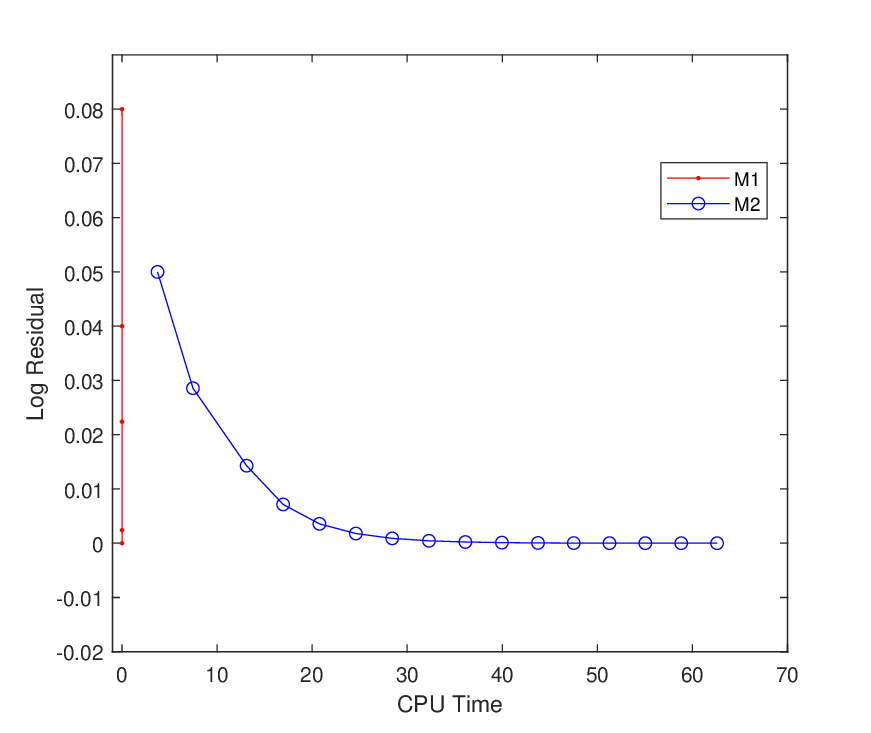}
\end{minipage}}
\caption{Errors and CPU times of Table 6.}\label{f4}
\end{figure}

In this example, $\mathrm{H}_{1}$ is symmetric positive definite. We
take $w=5$ for Method \ref{m32}. In such case, Tab. 6 reports the
numerical results for Example 5.5. 

\vspace{0.3cm}

From numerical results in Tabs. 5.5 and 5.6, the proposed Method 3.2
requires less iterative steps and CPU times than Method 3.3, which
implies that the former is significantly superior to the latter in
terms of computational efficiency.

\section{Conclusion}

This paper addresses two fundamental challenges in the study of the
EHLCP: the lack of iterative methods and error bounds. By leveraging
a max-min transformation technique, we derive an equivalent
fixed-point formulation that enables the development of efficient
iterative algorithms and global error analysis. The proposed methods
avoid pivoting operations, making them suitable for large-scale
sparse problems. Theoretical results are validated through numerical
experiments, demonstrating practical applicability. Future work will
focus on designing accelerated iterative schemes and exploring
applications in broader domains.


\section*{Funding}
This research is supported by National Natural Science Foundation of
China (No.12661083), Cross-integration Innovation team of modern
Applied Mathematics and Life Sciences in Yunnan Province, China
(202405AS350003), Yunnan Key Laboratory of Modern Analytical
Mathematics and Applications (202302AN360007).
\section*{Data availability} The data that support the findings
of this study are available from the corresponding author upon request.
\section*{Conflict of interest}
There is no conflict of interest.

\end{document}